\documentclass[a4paper, 11pt,reqno]{amsart} 

\usepackage{amsmath,amsfonts,amssymb,amsthm}
\usepackage[margin=1 in]{geometry}
\usepackage{esint} 

\usepackage{color}

\providecommand{\ve}{v^\epsilon}

\providecommand{\dx}{\, \mathrm{d} x}
\providecommand{\dy}{\, \mathrm{d} y}

\providecommand{\dt}{\, \mathrm{d} t}
\providecommand{\dr}{\, \mathrm{d} r}
\providecommand{\ds}{\, \mathrm{d} s}
\newcommand{\en}[1]{\left< #1 \right>}  

\newcommand{\norm}[1]{\lVert#1\rVert}
\newcommand{\supp}{\mathop{\mathrm{supp}}}

\renewcommand{\a}[1]{\left| #1 \right|}

\providecommand{\ah}{a_{\mathrm{hom}}}

\providecommand{\C}{\mathcal{C}}

\newtheorem{theorem}{Theorem}
\newtheorem{lemma}[theorem]{Lemma}

\newtheorem{proposition}[theorem]{Proposition}

\begin{document}

\author{Peter Bella}
\address{Universit\"at Leipzig, Mathematisches Institut}
\thanks{The first author is supported by the German Science Foundation DFG in context of the Emmy Noether junior research group BE 5922/1-1.}
\email{bella@math.uni-leipzig.de}

\author{Alberto Chiarini}
\address{Universit\'e d'Aix-Marseille}
\thanks{The second author is supported by LabEx Archim\'ede}
\email{alberto.chiarini@univ-amu.fr}

\author{Benjamin Fehrman}
\thanks{The third author is supported by the National Science Foundation Mathematical Sciences Postdoctoral Research Fellowship under Grant Number 1502731.}
\address{Max Planck Institute for Mathematics in the Sciences}
\email{fehrman@mis.mpg.de}

\title{A Liouville theorem for stationary and ergodic ensembles of parabolic systems}

\maketitle

\begin{abstract}  A first-order Liouville theorem is obtained for random ensembles of uniformly parabolic systems under the mere qualitative assumptions of stationarity and ergodicity.  Furthermore, the paper establishes, almost surely, an intrinsic large-scale $\C^{1,\alpha}$-regularity estimate for caloric functions.

\end{abstract}

\section{Introduction and main results.}  This paper considers random ensembles of uniformly parabolic systems
\begin{equation}\label{i_eq}u_t=\nabla\cdot a\nabla u,\end{equation}
where the law of the coefficient field $a$ is assumed to be \emph{stationary} with respect to space-time translations and \emph{ergodic}.  Precisely, for a probability space of coefficient fields $(\Omega,\mathcal{F},\en{\cdot})$, where $\en{\cdot}$ is used simultaneously to denote the law and expectation of the ensemble, the stationarity asserts that the coefficients are statistically homogeneous in time and space in the sense that
\begin{equation}\label{i_stationary}
\forall x\in\mathbb{R}^d, \forall t\in\mathbb{R}: \qquad a(\cdot, \cdot)\;\;\textrm{and}\;\;a(\cdot+x, \cdot+t)\;\;\textrm{have the same law under} \en{\cdot}.
\end{equation}
The ergodicity asserts that every translationally invariant function of the coefficient field is constant.  That is, for every bounded random variable $F$: 
\begin{equation}\label{i_ergodic}
\textrm{ if } \forall x\in\mathbb{R}^d, \forall t\in\mathbb{R}, \textrm{ and for } \en{\cdot}\textrm{-a.e.}\;a : F(a)=F(a(\cdot+x,\cdot+t)), \; \textrm{ then } F = c \; \en{\cdot}\textrm{-a.s.}\end{equation}
Finally, the ensemble is bounded and uniformly elliptic in the sense that there exists a deterministic $\lambda\in(0,1]$ such that 
\begin{equation}\label{i_bounded}\a{a\xi}\leq\a{\xi}\;\;\textrm{and}\;\;\lambda\a{\xi}^2\leq\xi\cdot a\xi \qquad \forall \xi \in \mathbb{R}^d, \textrm{ and for }\en{\cdot}\textrm{-a.e. } a.\end{equation}
Assumptions (\ref{i_stationary}) and (\ref{i_ergodic}) are the minimal statistical requirements on the ensemble $\en{\cdot}$ which guarantee the qualitative homogenization of equations like (\ref{i_eq}), see (\ref{i_qualitative_homogenization}).  Their role in this paper, and in homogenization theory generally, appears most essentially through applications of the ergodic theorem.  See, for instance, the foundational work of Papanicolaou and Varadhan \cite{papvar}, who worked in the elliptic setting.

However, conditions (\ref{i_stationary}) and (\ref{i_ergodic}) are merely qualitative and contain no quantitative information about the mixing properties of the ensemble.  Therefore, while the results of this paper apply to a very general class of environments, the corresponding homogenization may occur at an arbitrarily slow rate.  In order to obtain more quantitative statements, such as in the recent work Armstrong, Bordas and Mourrat \cite{ArmstrongBordasMourrat}, it would be necessary to quantify the ergodicity in the way, for example, of a spectral gap inequality or a finite-range of dependence.

The qualitative theory of homogenization for systems like (\ref{i_eq}) aims to characterize, for $\en{\cdot}$-a.e. $a$, the limiting behavior, as $\epsilon\rightarrow 0$, of solutions to the rescaled equation
\begin{equation}\label{i_eq_rescaled}\left\{\begin{array}{ll} u^\epsilon_t=\nabla\cdot a^\epsilon\nabla u^\epsilon & \textrm{in}\;\;\mathbb{R}^d\times(0,\infty) \\ u^\epsilon=u_0 & \textrm{on}\;\;\mathbb{R}^d\times\{0\},\end{array}\right.\end{equation}
where
$$a^\epsilon(\cdot,\cdot):=a\left(\frac{\cdot}{\epsilon}, \frac{\cdot}{\epsilon^2}\right)$$
is a parabolic rescaling of the coefficient field.   This is understood classically through the introduction of a space-time corrector $\phi=\{\phi_i\}_{i\in\{1,\ldots,d\}}$ satisfying, for each $i\in\{1,\ldots,d\}$,
\begin{equation}\label{i_corrector}\phi_{i,t}=\nabla\cdot a(\nabla\phi_i+e_i)\;\;\textrm{in}\;\;\mathbb{R}^{d+1}.\end{equation}
Then, in view of the linearity, for each $\xi\in\mathbb{R}^d$ the corresponding corrector $\phi_\xi$ is defined by the sum
\begin{equation}\label{i_corrector_xi}\phi_\xi:=\xi_i\phi_i,\end{equation}
where here, and throughout the paper, the notation employs Einstein's summation convention over repeated indices.

The gradient of the corrector $\nabla\phi$ is a random field which is stationary with finite energy.  That is, for each $x\in\mathbb{R}^d$, $t\in\mathbb{R}$ and $a\in\Omega$,
$$\nabla\phi(x,t;a)=\nabla\phi(0,0;a(\cdot+x,\cdot+t)),$$
and, for each $i\in\{1,\ldots,d\}$,
$$\en{\a{\nabla\phi_i}^2}<\infty.$$
These facts are used to prove the strict sublinearity of the large-scale $L^2$-averages of $\phi$ on parabolic cylinders.  Namely, for each $R>0$, let $B_R$ denote the ball of radius $R$ centered at the origin and let $\mathcal{C}_R$ denote the parabolic cylinder
$$\mathcal{C}_R:=B_R\times(-R^2,0].$$
The corrector satisfies, for $\en{\cdot}$-a.e. $a$, for each $i\in\{1,\ldots,d\}$,
\begin{equation}\label{i_corrector_sub}\lim_{R\rightarrow\infty}\frac{1}{R}\left(\fint_{\mathcal{C}_R}\a{\phi_i}^2\right)^\frac{1}{2}=0,\end{equation}
where here, and throughout the paper, the integration variables will be omitted unless there is a possibility of confusion.  This sublinearity is essentially equivalent to homogenization, see (\ref{i_qualitative_homogenization}) below, and is crucial for the arguments of this paper.

The corrector is used to identify the homogenized coefficient field $\ah$ as the expectation of the components of the flux according to the rules, for $i\in\{1,\ldots,d\}$,
\begin{equation}\label{i_homogenized_coefficients}\ah e_i:=\en{a(\nabla\phi_i+e_i)},\end{equation}
where the flux $q=\{q_i\}_{i\in\{1,\ldots,d\}}$ is defined, for each $i\in\{1,\ldots,d\}$, by
\begin{equation}\label{i_flux} q_i:=a(\nabla\phi_i+e_i).\end{equation}
It is a classical fact that the homogenized coefficient field $\ah$ is uniformly elliptic and bounded, as is shown in Lemma~\ref{i_lemma1}.  The solution of the corresponding constant-coefficient parabolic equation
\begin{equation}\label{intro_hom}\left\{\begin{array}{ll} v_t=\nabla\cdot \ah\nabla v & \textrm{in}\;\;\mathbb{R}^d\times(0,\infty) \\ v=u_0 & \textrm{on}\;\;\mathbb{R}^d\times\{0\},\end{array}\right.\end{equation}
then characterizes the limiting behavior, for $\en{\cdot}$-a.e. $a$ and as $\epsilon\rightarrow 0$, of the solutions to (\ref{i_eq_rescaled}).  Indeed, by obtaining an energy estimate for the error in the asymptotic expansion
$$u^\epsilon \simeq v+\epsilon\phi_i\left(\frac{\cdot}{\epsilon},\frac{\cdot}{\epsilon^2}\right)\partial_iv,$$
which relies upon the sublinearity (\ref{i_corrector_sub}), it follows that, for $\en{\cdot}$-a.e. $a$, for every $u_0\in L^2(\mathbb{R}^d)$ and $T>0$, as $\epsilon\rightarrow 0$,
\begin{equation}\label{i_qualitative_homogenization}u^\epsilon\rightarrow v\;\;\textrm{strongly in}\;\;L^2(\mathbb{R}^d\times[0,T]).\end{equation}
This almost sure convergence is the qualitative homogenization of the original ensemble.

Looking ahead, observe that the behavior of the solution $u^\epsilon$ to (\ref{i_eq_rescaled}) on a unit scale, for $\epsilon>0$ small, corresponds to a characterization of the large-scale behavior of the solution $u$ satisfying (\ref{i_eq}).  Namely, the behavior of the solution $u^\epsilon$ on a unit scale corresponds to the behavior of $u$ on scale $\epsilon^{-1}$ in space and $\epsilon^{-2}$ in time.  The purpose of this paper will be to characterize the extent to which solutions of (\ref{i_eq}) inherit, on large-scales and for $\en{\cdot}$-a.e. $a$, the regularity of solutions to constant-coefficient parabolic equations.

A concise statement of this large-scale regularity is contained in the following first-order Liouville theorem, which is the main theorem of the paper.

\begin{theorem}\label{i_liouville}  Suppose that $\en{\cdot}$ is stationary (\ref{i_stationary}), ergodic (\ref{i_ergodic}) and bounded and uniformly elliptic (\ref{i_bounded}).  Then, $\en{\cdot}$-a.e. $a$ satisfies the following first-order Liouville property: if $u$ is an ancient whole-space $a$-caloric function, that is if $u$ is a distributional solution of
$$u_t=\nabla\cdot a\nabla u\;\;\textrm{in}\;\;\mathbb{R}^d\times(-\infty,0),$$
which is strictly subquadratic on parabolic cylinders in the sense that, for some $\alpha\in(0,1)$,
$$\lim_{R\rightarrow\infty}\frac{1}{R^{1+\alpha}}\left(\fint_{\mathcal{C}_R}\a{u}^2\right)^{\frac{1}{2}}=0,$$
then there exists $c\in\mathbb{R}$ and $\xi\in\mathbb{R}^d$ such that
$$u(x,t)=c+x\cdot\xi+\phi_\xi(x,t)\;\;\textrm{in}\;\;\mathbb{R}^d\times(-\infty,0),$$
for the corrector $\phi_\xi$ defined in (\ref{i_corrector_xi}). 
\end{theorem}

The proof of Theorem~\ref{i_liouville} is strongly motivated by the work of Gloria, Neukamm and Otto \cite{GNO4}, who considered precisely these questions for stationary and ergodic ensembles of elliptic equations.  It is based on controlling the large-scale $L^2$-deviation of the gradient of an $a$-caloric function from the span of the $a$-caloric gradients $\{\xi+\nabla \phi_\xi\}_{\xi\in\mathbb{R}^d}$.  The excess of an $a$-caloric function measures this deviation, and is defined, for each $R>0$ and $a$-caloric function $u$ on $\mathcal{C}_R$, by
\begin{equation}\label{i_excess}\textrm{Exc}(u;R):=\inf_{\xi\in\mathbb{R}^d}\fint_{\mathcal{C}_R}(\nabla u-\xi-\nabla\phi_\xi)\cdot a(\nabla u-\xi-\nabla\phi_\xi).\end{equation}
In Proposition~\ref{i_excess_decay} below, for $\en{\cdot}$-a.e. $a$, the excess of an $a$-caloric function will be shown to decay like a power law in the radius.  However, before the statement, it is useful to observe some essential differences between the parabolic and elliptic settings.  In what follows, the superscript ``$\textrm{ell}$'' will be used to differentiate elliptic objects from their parabolic counterparts.

In the elliptic case, for a stationary and ergodic ensemble $\en{\cdot}^\textrm{ell}$ of bounded, uniformly elliptic coefficient fields $a^\textrm{ell}$, the corrector $\phi^\textrm{ell}=\{\phi^{\textrm{ell}}_i\}_{i\in\{1,\ldots,d\}}$ is defined by the equations, for $i\in\{1,\ldots,d\}$,
\begin{equation}\label{i_elliptic_corrector}-\nabla\cdot a^\textrm{ell}(\nabla\phi^\textrm{ell}_i+e_i)=0\;\;\textrm{in}\;\;\mathbb{R}^d,\end{equation}
and, for each $\xi\in\mathbb{R}^d$,
$$\phi^\textrm{ell}_\xi:=\phi^\textrm{ell}_i\xi_i.$$
These correctors play a virtually identical role to the parabolic correctors (\ref{i_corrector}) in elliptic homogenization theory.

Excess for uniformly elliptic ensembles was first defined in \cite[Lemma~2]{GNO4}, although they worked not with the intrinsic energy defined by $a$ but with the equivalent $L^2$-energy.  This differs, for instance, from the definition used in the work of the first author, third author and Otto \cite{BellaFehrmanOtto}, which considered degenerate elliptic ensembles for which it was essential to incorporate the environment $a$.  These notions motivated definition (\ref{i_excess}), and measured the deviation of the gradient of an $a^\textrm{ell}$-harmonic function $u$ on $B_R$, by which is meant a solution
$$-\nabla\cdot a^\textrm{ell}\nabla u=0\;\;\textrm{in}\;\;B_R,$$
from the span of $a^{\textrm{ell}}$-harmonic gradients $\{\xi+\nabla\phi^\textrm{ell}_\xi\}_{\xi\in\mathbb{R}^d}$.  Precisely, for each $R>0$ and $a^\textrm{ell}$-harmonic function $u$ on $B_R$,
\begin{equation}\label{i_elliptic_excess}\textrm{Exc}^\textrm{ell}(u;R):=\inf_{\xi\in\mathbb{R}^d}\fint_{B_R}(\nabla u-\xi-\nabla\phi^\textrm{ell}_\xi)\cdot a^\textrm{ell}(\nabla u-\xi-\nabla\phi^\textrm{ell}_\xi).\end{equation}

The decay of the excess was controlled in \cite[Lemma~2]{GNO4}  through the introduction of a flux correction $\sigma^\textrm{ell}=\{\sigma^\textrm{ell}_i\}_{i\in\{1,\ldots,d\}}$.  Namely, the flux $q^\textrm{ell}=\{q^\textrm{ell}_i\}_{i\in\{1,\ldots,d\}}$ is defined, for each $i\in\{1,\ldots,d\}$, by
$$q^\textrm{ell}_i:=a^\textrm{ell}(\nabla\phi^\textrm{ell}_i+e_i),$$
where, in analogy with the parabolic setting, the homogenized coefficient field $\ah^\textrm{ell}$ is defined by the expectation of the components of the flux, for $i\in\{1,\ldots,d\}$,
$$\ah^\textrm{ell} e_i:=\en{a^\textrm{ell}(\nabla\phi^\textrm{ell}_i+e_i)}^\textrm{ell}.$$
Therefore, strictly in the elliptic case, the corrector equation (\ref{i_elliptic_corrector}) asserts that the components of the flux are divergence-free and may be viewed as closed $(d-1)$-forms on the whole space.  Hence, for each $i\in\{1,\ldots,d\}$, there exists a $(d-2)$-form, which is represented by a skew-symmetric matrix $\sigma^\textrm{ell}_i=(\sigma_{ijk})^\textrm{ell}_{j,k\in\{1,\ldots,d\}}$, satisfying
\begin{equation}\label{intro_sigma} \nabla\cdot\sigma^\textrm{ell}_i=q^{\textrm{ell}}_i-\ah^\textrm{ell} e_i,\end{equation}
where the divergence of the tensor-field $\sigma_i^\textrm{ell}$ is defined, for each $i,j\in\{1,\ldots,d\}$, by
$$(\nabla\cdot \sigma^\textrm{ell}_i)_j=\sum_{k=1}^d\partial_k\sigma^\textrm{ell}_{ijk}.$$
Furthermore, the flux correction $\sigma^\textrm{ell}$ is fixed according to the choice of gauge, for each $i,j,k\in\{1,\ldots,d\}$,
$$\Delta \sigma^\textrm{ell}_{ijk}=\partial_kq^\textrm{ell}_{ij}-\partial_jq^\textrm{ell}_{ik}.$$
In \cite[Lemma~2]{GNO4}, the sublinearity of the large-scale $L^2$-averages of the extended corrector $(\phi^\textrm{ell},\sigma^\textrm{ell})$ on large balls is shown to imply that the elliptic excess (\ref{i_elliptic_excess}) decays as a power law in the radius.

Precisely, for each $\alpha\in(0,1)$, there exists $C^\textrm{ell}_0=C^\textrm{ell}_0(\alpha,d,\lambda)>0$ and $C_1^\textrm{ell}=C_1^\textrm{ell}(\alpha, d,\lambda)>0$ for which, whenever a pair of radii $0<r<R<\infty$ satisfy, for every $\rho\in[r,R]$, for each $i\in\{1,\ldots,d\}$,
$$\frac{1}{\rho}\left(\fint_{B_\rho}\a{\phi_i^\textrm{ell}}^2+\a{\sigma^\textrm{ell}_i}^2\right)^{\frac{1}{2}}\leq\frac{1}{C^\textrm{ell}_0},$$
then, for every $a^\textrm{ell}$-harmonic function $u$ in $B_R$,
\begin{equation}\label{i_elliptic_decay}\textrm{Exc}^\textrm{ell}(u;r)\leq C^\textrm{ell}_1\left(\frac{r}{R}\right)^{2\alpha}\textrm{Exc}^\textrm{ell}(u;R).\end{equation}
Observe, in particular, that this is a deterministic result.  Indeed, the stochastic properties of the extended corrector $(\phi^\textrm{ell},\sigma^\textrm{ell})$ are necessary to prove that, for $\en{\cdot}$-a.e. $a$, the large-scale $L^2$-averages are sublinear.  But, by taking this fact as an input, it follows from a Campanato iteration that the excess of an arbitrary $a^\textrm{ell}$-harmonic function decays according to (\ref{i_elliptic_decay}).  In this paper, the analogous result will also be obtained for the parabolic excess, as shown in Proposition~\ref{i_excess_decay} below.

The first essential difference is that, unlike in the elliptic case, the fluxes $\{q_i\}_{i\in\{1,\ldots,d\}}$ defined in (\ref{i_flux}) are not divergence-free, and so an immediate analogue of the flux correction $\sigma^\textrm{ell}=\{\sigma^\textrm{ell}_i\}_{i\in\{1,\ldots,d\}}$ cannot be defined.  Instead, the flux is essentially decomposed according to the Weyl decomposition, where the parabolic $\sigma$ is constructed to correct the divergence-free component.  Precisely, for each $i\in\{1,\ldots,d\}$,
$$q_i=q_{i,\textrm{sol}}+\nabla\psi_i,$$
where the solenoidal part $q_{i,\textrm{sol}}$ is divergence-free and $\nabla\psi_i$ is a stationary, finite-energy gradient chosen to satisfy
\begin{equation}\label{i_psi}\Delta\psi_i=\nabla\cdot q_i.\end{equation}
Indeed, for each $i\in\{1,\ldots,d\}$, one first defines $\nabla\psi_i$ according to (\ref{i_psi}) and then observes that $q_i-\nabla\psi_i$ is divergence-free.

The flux correction $\sigma=\{\sigma_i\}_{i\in\{1,\ldots,d\}}$ is then defined, for each $i\in\{1,\ldots,d\}$, by the equation
\begin{equation}\label{i_sigma}\nabla\cdot\sigma_i=\left(q_i-\nabla\psi_i\right)-\en{q_i\;|\;\mathcal{F}_{\mathbb{R}^d}},\end{equation}
where $\en{q_i\;|\;\mathcal{F}_{\mathbb{R}^d}}$ denotes the conditional expectation of $q_i$ with respect to the sub-sigma-algebra $\mathcal{F}_{\mathbb{R}^d}\subset\mathcal{F}$ of subsets of $\Omega$ which are invariant with respect to spatial translations of the coefficient fields.  They are fixed following the choice of gauge, for each $i,j,k\in\{1,\ldots,d\}$,
\begin{equation}\label{i_gauge}\Delta\sigma_{ijk}=\partial_k(q_i-\nabla\psi_i)_j-\partial_j(q_i-\nabla\psi_i)_k.\end{equation}

Finally, for each $i\in\{1,\ldots,d\}$, it is necessary to correct the oscillations of the conditional expectation $\en{q_i\;|\;\mathcal{F}_{\mathbb{R}^d}}$ about its mean.  The corrector $\zeta=\{\zeta_{ij}\}_{i,j\in\{1,\ldots,d\}}$ is constructed explicitly for this purpose and satisfies, for each $i\in\{1,\ldots,d\}$,
\begin{equation}\label{i_zeta} \partial_t \zeta_i=\en{q_i\;|\;\mathcal{F}_{\mathbb{R}^d}}-\en{q_i}=\en{q_i\;|\;\mathcal{F}_{\mathbb{R}^d}}-\ah e_i.\end{equation}
In particular, this final correction $\zeta$ is constant in space, as a $\mathcal{F}_{\mathbb{R}^d}$-measurable field,  and depends only on time.

In comparison with the elliptic setting, where the decay of the excess was determined by the sublinearity of the large-scale $L^2$-averages of $(\phi^\textrm{ell},\sigma^\textrm{ell})$, the decay of the parabolic excess will be determined by the sublinearity of the large-scale $L^2$-averages of the corrector $(\phi,\psi,\sigma)$, measured with respect to the scaling in space, and the sublinearity of the large-scale $L^2$-averages of $\zeta$, measured with respect to the scaling in time.  The first lemma of the paper establishes the existence of the extended corrector $(\phi,\psi,\sigma,\zeta)$.

\begin{lemma}\label{i_lemma1} Suppose that the ensemble $\en{\cdot}$ satisfies (\ref{i_stationary}), (\ref{i_ergodic}) and (\ref{i_bounded}).  There exist $C=C(d,\lambda)>0$ and four random fields $\phi=\{ \phi_i \}_{i\in\{1,\ldots,d\}}$, $\psi=\{\psi_i\}_{i\in\{1,\ldots,d\}}$, $\sigma=\{\sigma_{ijk}\}_{i,j,k\in\{1,\ldots,d\}}$ and $\zeta=\{\zeta_{ij}\}_{i,j\in\{1,\ldots,d\}}$ on $\mathbb{R}^{d+1}$ with the following properties:

The gradient fields are stationary, finite energy random processes with vanishing expectation:  for each $i,j,k\in\{1,\ldots,d\}$,
$$\en{ \a{\nabla\phi_i}^2 } + \en{ \a{\nabla \psi_i}^2} + \en{ \a{\nabla \sigma_{ijk}}^2} +\en{\a{\partial_t\zeta_{ij}}^2} \leq C,$$
and
$$\en{\nabla \phi_i} = \en{\nabla\psi_i} = \en{\nabla \sigma_{ijk} } =\en{\partial_t \zeta_{ij}} = 0.$$
For each $i\in\{1,\ldots,d\}$, the field $\sigma_i=(\sigma_{ijk})_{j,k\in\{1,\ldots,d\}}$ is skew-symmetric in its last two indices:  for each $i,j,k\in\{1,\ldots,d\}$,
$$\sigma_{ijk} = -\sigma_{ikj}.$$
The fields $\psi$ and $\sigma$ are stationary in time:  for each $x\in\mathbb{R}^d$, $t\in\mathbb{R}$ and $a\in\Omega$,
$$\psi(x,t;a)=\psi(x,0;a(\cdot,\cdot+t)),$$
and
$$\sigma(x,t;a)=\sigma(x,0;a(\cdot,\cdot+t)).$$

Furthermore, for $\en{\cdot}$-a.e. $a$, the following equations are satisfied in the sense of distributions on $\mathbb{R}^{d+1}$.  The field $\phi$ satisfies (\ref{i_corrector}):  for each $i\in\{1,\ldots,d\}$,
$$\phi_{i,t}=\nabla \cdot a(\nabla\phi_i+e_i).$$

The potential part of the flux is corrected by $\psi$ according to (\ref{i_psi}):  for each $i\in\{1,\ldots,d\}$,
$$\Delta\psi_i=\nabla\cdot q_i.$$

The field $\sigma$ corrects the divergence-free part of the flux according to (\ref{i_sigma}):  for each $i\in\{1,\ldots,d\}$,
$$ \nabla \cdot \sigma_i= q_i -\nabla\psi_i-\en{q_i\;|\;\mathcal{F}_{\mathbb{R}^d}},$$
where $\en{\cdot\;|\;\mathcal{F}_{\mathbb{R}^d}}$ denotes the conditional expectation with respect to the sub-sigma-algebra $\mathcal{F}_{\mathbb{R}^d}\subset\mathcal{F}$ of subsets of $\Omega$ which are invariant with respect to spatial translations of the coefficient field.  Furthermore, $\sigma$ is constructed according to the choice of gauge, for each $i,j,k\in\{1,\ldots,d\}$,
$$\Delta\sigma_{ijk}=\partial_k(q_i-\nabla\psi_i)_j-\partial_j(q_i-\nabla\psi_i)_k.$$

The field $\zeta$ corrects the oscillation of the conditional expectation about its mean:  for each $i\in\{1,\ldots,d\}$, the random vector field $\zeta_i$ is constant in space and satisfies
$$\partial_t\zeta_i=\en{q_i\;|\;\mathcal{F}_{\mathbb{R}^d}}-\en{q_i}=\en{q_i\;|\;\mathcal{F}_{\mathbb{R}^d}}-\ah e_i.$$

Finally, the homogenized coefficient field $\ah$ defined in (\ref{i_homogenized_coefficients}) is bounded and uniformly elliptic:  for each $\xi\in\mathbb{R}^d$,
$$\lambda \a{\xi}^2\leq \xi\cdot \ah \xi\;\;\textrm{and}\;\;\a {\ah \xi}\leq \frac{1}{\lambda} \a{\xi}.$$
\end{lemma}

The following two propositions effectively split the probabilistic and deterministic aspects of the paper.  Proposition \ref{i_sublinear} contains the probabilistic parts, and uses the stationarity and ergodicity of the ensemble to prove that the large-scale $L^2$-averages of $(\phi,\psi,\sigma)$ are sublinear with respect to the spatial scaling and that those of $\zeta$ are sublinear with respect to the time scaling.  This fact is essentially classical for the case of the correctors $\phi$ and $\zeta$, although a new argument for the sublinearity of $\phi$ is presented which may be of independent interest.  A new argument is required to prove the sublinearity of $\sigma$ and $\psi$.

The difference is the following.  The corrector $\phi$ is, in general, not stationary in either space or time but equation (\ref{i_corrector}) yields some control over both its spatial and temporal derivatives.  Similarly, the corrector $\zeta$ has an explicit, stationary time derivative but is itself not stationary.  In the second case, since equations (\ref{i_psi}) and (\ref{i_sigma}) yield only the spatial regularity for $\psi$ and $\sigma$, it is necessary to use the fact that both fields are stationary in time in order to obtain the convergence.

In fact, the following proposition will prove the sublinearity of the normalized corrector where, in the case of $\phi$, the components are normalized by their large-scale averages on a parabolic cylinder, in the case of $\zeta$, using the fact that the Sobolev embedding implies that $\zeta$ is continuous,  the components are normalized by their value at time zero and, in the case of $\psi$ and $\sigma$, the functions are normalized,  for each fixed time, by their large-scale averages on a ball.  This is in fact equivalent to the sublinearity of the corrector $(\phi,\psi,\sigma,\zeta)$ without a normalization, see for instance \cite[Lemma~2]{BellaFehrmanOtto}, but since this observation is not necessary for the arguments of the paper it is omitted.

For an arbitrary function $\varphi:\mathbb{R}^{d+1}\rightarrow\mathbb{R}$, define, for each $R>0$ and $t\in\mathbb{R}$,
$$(\varphi)_R:=\fint_{\mathcal{C}_R}\varphi\;\;\textrm{and}\;\;(\varphi)_{t,R}:=\fint_{B_R}\varphi(\cdot,t).$$
The precise normalization considered and the corresponding sublinearity are contained in the following proposition.

\begin{proposition}\label{i_sublinear} Suppose that the ensemble $\en{\cdot}$ satisfies (\ref{i_stationary}), (\ref{i_ergodic}) and (\ref{i_bounded}).  Then, for $\en{\cdot}$-a.e $a$, the corrector $(\phi,\psi,\sigma)$ is strictly sublinear with respect to the spatial scaling and the corrector $\zeta$ is strictly sublinear with respect to the time scaling in the sense that, for each $i,j,k\in\{1,\ldots,d\}$,
\begin{equation}\label{i_sublinear_1}\lim_{R\rightarrow\infty}\frac{1}{R}\left(\fint_{C_R}\a{\phi_i-(\phi_i)_R}^2+\a{\psi_i-(\psi_i)_{t,R}}^2+\a{\sigma_{ijk}-(\sigma_{ijk})_{t,R}}^2\right)^\frac{1}{2}=0, \end{equation}
and
\begin{equation}\label{i_sublinear_10}\lim_{R\rightarrow\infty}\frac{1}{R^2}\left(\fint_{C_R}\a{\zeta_{ij}-\zeta_{ij}(0)}^2\right)^\frac{1}{2}=0.\end{equation}
Furthermore, for $\en{\cdot}$-a.e. $a$, for each $i\in\{1,\ldots,d\}$, the large-scale $L^2$-averages of the components of the flux satisfy
\begin{equation}\label{i_sublinear_2}\lim_{R\rightarrow\infty}\left(\fint_{\mathcal{C}_R}\a{q_i}^2\right)^\frac{1}{2}=\en{\a{q_i}^2}^\frac{1}{2}.\end{equation}
\end{proposition}

It is important to observe at this point that equations (\ref{i_psi}) and (\ref{i_sigma}) defining $\psi$ and $\sigma$ are invariant if either $\psi$ or $\sigma$ is altered by a time stationary constant.  This explains why in (\ref{i_sublinear_1}), for each $R>0$, it remains useful to allow for a time-dependent normalization.  The equations (\ref{i_corrector}) and (\ref{i_zeta}) defining $\phi$ and $\zeta$ are not likewise invariant, and therefore the corresponding normalizations appearing in (\ref{i_sublinear_1}) and (\ref{i_sublinear_10}) are necessarily achieved by subtracting a true constant.

The deterministic aspect of the paper uses a Campanato iteration, which takes the conclusion of Proposition~\ref{i_sublinear} as input.  Namely, it will be shown that, for any $\alpha\in(0,1)$, the parabolic excess decays like a power law in the radius as soon as the quantities appearing in (\ref{i_sublinear_1}) and (\ref{i_sublinear_10}) are sufficiently small and as soon as (\ref{i_sublinear_2}) is sufficiently close to its expectation.  This is to say that there exists a random but $\en{\cdot}$-a.e. $a$ finite radius $r_*(a)$ such that, whenever $r_*<r<R<\infty$, for every $a$-caloric function $u$ in $\mathcal{C}_R$, the parabolic excess satisfies, for $C_1=C_1(\alpha,d,\lambda)>0$,
$$\textrm{Exc}(u;r)\leq C_1\left(\frac{r}{R}\right)^{2\alpha}\textrm{Exc}(u;R).$$
This is the content of the following proposition.

\begin{proposition}\label{i_excess_decay}  Suppose that the ensemble $\en{\cdot}$ satisfies (\ref{i_stationary}), (\ref{i_ergodic}) and (\ref{i_bounded}).  Fix a H\"older exponent $\alpha\in(0,1)$.  Then, there exist constants $C_0=C_0(\alpha, d,\lambda)>0$ and $C_1(\alpha,d,\lambda)>0$ with the following property:

If $R_1<R_2$ are two radii such that, for each $R\in[R_1,R_2]$ and for each $i,j,k\in\{1,\ldots,d\}$,
$$\frac{1}{R}\left(\fint_{\mathcal{C}_R}\a{\phi_i-(\phi_i)_R}^2+\a{\psi_i-(\psi_i)_{t,R}}^2+\a{\sigma_{ijk}-(\sigma_{ijk})_{t,R}}^2\right)^{\frac{1}{2}}\leq\frac{1}{C_0},$$
and
$$\frac{1}{R^2}\left(\fint_{\mathcal{C}_R}\a{\zeta_{ij}-\zeta_{ij}(0)}^2\right)^{\frac{1}{2}}\leq\frac{1}{C_0},$$
and such that, for each $i\in\{1,\ldots,d\}$ and $R\in[R_1,R_2]$,
$$\left(\fint_{\mathcal{C}_R}\a{q_i}^2\right)^{\frac{1}{2}}\leq 2\en{\a{q_i}^2}^\frac{1}{2},$$
then any distributional solution $u$ to the parabolic equation
$$u_t=\nabla\cdot a\nabla u\;\;\textrm{in}\;\;\mathcal{C}_{R_2}$$
satisfies
$$\mathrm{Exc}(u;R_1)\leq C_1\left(\frac{R_1}{R_2}\right)^{2\alpha} \mathrm{Exc}(u;R_2).$$

\end{proposition}

The proof of Proposition~\ref{i_excess_decay} is motivated by the proof of \cite[Lemma~2]{GNO4} from the elliptic setting.  There, the flux correction $\sigma^\textrm{ell}$ was used to express the residuum of the homogenization error in a useful divergence form.  That is, for $R>0$, given an $a^\textrm{ell}$-harmonic function $u$ in $B_R$, define its $\ah^\textrm{ell}$-harmonic extension $v$ into $B_R$ to be the solution
$$\left\{\begin{array}{rll} -\nabla\cdot \ah^{\textrm{ell}}\nabla v & =0 & \textrm{in}\;\;B_R \\ v & =u & \textrm{on}\;\;\partial B_R.\end{array}\right.$$
Then, for a smooth cutoff function $\eta$ vanishing along the  boundary $\partial B_R$, define the augmented homogenization error $w^\textrm{ell}$ to be the following modification of the classical two-scale expansion
$$w^\textrm{ell}:=u-(1+\eta \phi^\textrm{ell}_i\partial_i)v,$$
where the cutoff is used in order to guarantee the difference $w^\textrm{ell}$ vanishes on the boundary.

It was proven that the augmented homogenization error $w^\textrm{ell}$ satisfies
\begin{equation}\label{i_elliptic_residuum}\left\{\begin{array}{rll} -\nabla\cdot a^\textrm{ell}\nabla w^\textrm{ell}=& \nabla\cdot \left((1-\eta)(a^\textrm{ell}-\ah^\textrm{ell})\nabla v+(\phi^\textrm{ell}_ia^\textrm{ell}-\sigma^\textrm{ell}_i)\nabla(\eta\partial_iv)\right) & \textrm{in}\;\;B_R \\ w^\textrm{ell}= & 0 & \textrm{on}\;\;\partial B_R,\end{array}\right.\end{equation}
which, by testing the equation with $w^\textrm{ell}$,  yields a useful energy estimate that provides the starting point for a Campanato iteration.

In particular, by analyzing the right hand side of (\ref{i_elliptic_residuum}), the energy of $w^\textrm{ell}$ can be controlled by the growth of the extended corrector $(\phi^\textrm{ell},\sigma^\textrm{ell})$, the choice of the cutoff function $\eta$ and the interior and boundary regularity of the $\ah^\textrm{ell}$-harmonic function $v$.  The argument is completed by observing that, owing to the regularity of $\ah^\textrm{ell}$-harmonic functions, the energy of the homogenization error is a good approximation for the excess.

The methods of this paper apply the same philosophy to the parabolic setting.  However, similarly to what was done in the proof of \cite[Theorem~2]{BellaFehrmanOtto}, it is furthermore necessary to introduce a spatial regularization of the $a$-caloric function $u$.  The purpose of this is to obtain more regularity in time, since such functions are already sufficiently regular in space.  Precisely, if $u$ is an $a$-caloric function then, in general, its time derivative $u_t\in H^{-1}$ and no better, where $H^{-1}$ denotes the dual space of the Sobolev space $H^1$.  However, for every $\epsilon>0$, if $u^\epsilon$ denotes the spatial convolution of $u$ on scale $\epsilon>0$, then it is possible to show that $u^\epsilon_t\in L^2$, see Section~5.1 below.  This additional approximation is necessary in order to apply the boundary estimate of Section~5.3.

For $\epsilon>0$, the $\ah$-caloric function $v^\epsilon$ will then be the $\ah$-caloric extension of the spatial regularization $u^\epsilon$ into $\mathcal{C}_R$.  Namely, for $R>0$ and $\epsilon>0$, given an $a$-caloric function $u$ in $\mathcal{C}_{R+\epsilon}$, define the $\ah$-caloric extension $v^\epsilon$ of $u^\epsilon$ into $\mathcal{C}_R$ to be the solution
$$\left\{\begin{array}{rll} \ve_t & =\nabla\cdot\ah\nabla \ve & \textrm{in}\;\;\mathcal{C}_R \\ \ve & =u^\epsilon & \textrm{on}\;\;\partial_p\mathcal{C}_R,\end{array}\right.$$
where $\partial_p\mathcal{C}_R$ denotes the parabolic boundary
$$\partial_p\mathcal{C}_R:=\left(B_R\times\{-R^2\}\right)\cup\left(\partial B_R\times[-R^2,0]\right).$$
Then, again motivated by the classical two-scale expansion, for $\epsilon>0$ and a smooth cutoff function $\eta$ vanishing on the parabolic boundary $\partial_p\mathcal{C}_R$, the augmented homogenization error $w$ will be defined as
$$w:=u-(1+\eta\phi_i\partial_i)v^\epsilon.$$
It will be shown in the proof of Proposition~\ref{i_excess_decay} that the augmented homogenization error satisfies
\begin{equation}\label{i_residuum}\left\{\begin{array}{rll} w_t-\nabla\cdot a\nabla w  = & \nabla\cdot\left((1-\eta)(a-\ah)\nabla \ve\right)+\nabla\cdot\left((\phi_ia+\psi_i-\sigma_i)\nabla(\eta\partial_i\ve)\right) & \textrm{in}\;\;\mathcal{C}_R \\    & + \partial_t\zeta_i\cdot\nabla(\eta\partial_i\ve)-\phi_i (\eta\partial_i\ve)_t-\psi_i\Delta(\eta\partial_i\ve) & \\ w  = & u-u^\epsilon &   \textrm{on}\;\;\partial_p\mathcal{C}_R.\end{array}\right.\end{equation}
As in the elliptic setting, the energy estimate obtained by testing this equation with $w$, for an appropriately chosen cutoff $\eta$, will be the starting point of the Campanato iteration used to control the decay of the excess.  In this case, there is a contribution from the boundary, which will be controlled first by fixing $\epsilon>0$ small.  From the right hand side of (\ref{i_residuum}), the energy of the homogenization error will then be controlled by the growth of the extended parabolic corrector $(\phi,\psi,\sigma)$ and, after integrating in parts by time, the growth of $\zeta$ and $q$.  It is furthermore necessary to make a good choice for the cutoff function $\eta$ and to use the interior and boundary regularity of the $\ah$-caloric function $\ve$.  The argument is completed by observing that, owing to the interior regularity of $\ah$-caloric functions, the homogenization error $w$ provides a good approximation for the excess.

Finally, the following parabolic Caccioppoli inequality will be used in the proofs of Theorem~\ref{i_liouville} and Proposition~\ref{i_excess_decay}.  The proof is classical, and is included for the convenience of the reader.

\begin{lemma}\label{i_cacc}  Suppose that $\en{\cdot}$ satisfies (\ref{i_bounded}).  There exists $C=C(\lambda)>0$ such that, for $\en{\cdot}$-a.e. $a$, for every $R>0$ and distributional solution $u$ of the equation
$$u_t=\nabla\cdot a\nabla u\;\;\textrm{in}\;\;\mathcal{C}_R,$$
and for every $c\in\mathbb{R}$ and $\rho\in(0,\frac{R}{2})$,
\begin{equation}\label{i_cacc_eq}\int_{\mathcal{C}_{R-\rho}}\a{\nabla u}^2\leq \frac{C}{\rho^2}\int_{\mathcal{C}_R\setminus \mathcal{C}_{R-\rho}}\a{u-c}^2.\end{equation}
\end{lemma}

In comparison with the elliptic setting, the qualitative homogenization theory of divergence-form operators with coefficients depending on time and space is relatively under studied.  While the case of periodic coefficients has long been understood, and a full explanation can be found in the classic reference Bensoussan, Lions and Papanicolaou \cite[Chapter~3]{BensoussanLionsPapanicolaou}, the qualitative stochastic homogenization of stationary and ergodic ensembles like (\ref{i_eq}) was obtained only more recently by Rhodes \cite{Rhodes1,Rhodes2}.  However, related problems were earlier handled, such as the case of a Brownian motion in the presence of a divergence-free drift, by Komorowski and Olla \cite{KomorowskiOlla}, Landim, Olla and Yau \cite{LandimOllaYau} and Oelschl\"ager \cite{Oelschlager}.  In the discrete setting, related questions have been considered, for instance, by Andres \cite{Andres}, Bandyopadhyay and Zeitouni \cite{BandyopadhyayZeitouni} and Rassoul-Agha and Sepp\"al\"ainen \cite{Rassoul-AghaSepp} in the uniformly elliptic setting and, for degenerate environments, by Andres, the second author, Deuschel and Slowik \cite{ACDS}.

The quantitative homogenization of such ensembles has only recently been considered, and the preprint \cite{ArmstrongBordasMourrat} contains, to our knowledge, the first results in this direction.  In particular, in \cite[Theorem~1.2]{ArmstrongBordasMourrat}, a full hierarchy of Liouville theorems is obtained for ensembles satisfying a finite-range dependence in space and time.  Their method is motivated by the work of Armstrong and Smart \cite{armstrongsmart2014} from the elliptic setting, which adapted the approach of Avellaneda and Lin \cite{AvellanedaLinCPAM87} from the context of periodic homogenization.

In \cite{AvellanedaLinCPAM87}, a full hierarchy of Liouville properties was established for uniformly elliptic and periodic coefficient fields based upon the previous works Avellaneda and Lin \cite{AvellanedaLin87,AvellanedaLinCRAS89}, which developed a large-scale regularity theory in H\"older and $L^p$-spaces.  In \cite{armstrongsmart2014}, the approach of \cite{AvellanedaLinCPAM87} was adapted to stationary and ergodic ensembles satisfying a finite-range dependence.  Their proof, which obtained a large-scale $C^{0,1}$-regularity theory, was based upon a variational approach and the quantification of the convergence of certain sub-additive and super-additive energies.  Their work was later extended by Armstrong and Mourrat \cite{ArmstrongMourrat} to more general mixing conditions, and subsequently gave rise to a significant literature on the subject.  The interested reader is pointed to the recent monograph Armstrong, Kuusi and Mourrat \cite{ArmstrongKuusiMourrat}, and the references therein.

The approach of this paper follows closely the work \cite{GNO4}, which derived, for uniformly elliptic ensembles, a large-scale $C^{1,\alpha}$-regularity estimate and first-order Liouville property under the qualitative assumptions of stationarity and ergodicity.   The method was based upon the introduction of an intrinsic notion of excess, as defined in (\ref{i_elliptic_excess}), as well as the construction of the flux correction $\sigma^\textrm{ell}$ defined in (\ref{intro_sigma}).  The introduction of $\sigma^\textrm{ell}$ was used to prove that the homogenization error solves the  divergence-form equation (\ref{i_elliptic_residuum}), which provided the starting point for a Campanato iteration as explained above.

Subsequently, Fischer and Otto \cite{FischerOtto} obtained a  full hierarchy of Liouville properties under a mild quantification of the ergodicity.  In Fischer and Otto \cite{FischerOtto2}, the necessary quantification of ergodicity from \cite{FischerOtto} was shown to be satisfied by a general class of Gaussian environments.  However, absent some mild quantification of ergodicity in the sense of either \cite{ArmstrongMourrat} or \cite{FischerOtto}, the existence of higher order Liouville and large-scale regularity statements remains an open question.  

Finally, motivated by the work of the second author and Deuschel \cite{ChiariniDeuschel}, the first author, third author and Otto \cite{BellaFehrmanOtto} derived a large-scale $C^{1,\alpha}$-regularity theory and first-order Liouville theorem for degenerate elliptic equations, where the boundedness and uniform ellipticity (\ref{i_bounded}) was replaced by certain moment conditions.  It is expected that the results of this paper can be similarly extended to degenerate environments, and the setting of \cite{ACDS} will serve as the starting point for future work.

In principle, one could also hope to combine the methods of this paper with those of \cite{GNO4}, in the presence of a logarithmic Sobolev inequality like that used in \cite[Theorem~1]{GNO4}, to obtain more quantitative information.  For example, the minimal radius $r_*(a)>0$ defined, for $C_0>0$ from Proposition~\ref{i_excess_decay}, after suppressing the normalizations in the notation, by
\begin{multline*} r_*(a)=\inf\left\{\;r\geq1\;|\;\textrm{For all}\;\;R\geq r,\;\textrm{for each}\;i\in\{1,\ldots,d\},\;\;\left(\fint_{\mathcal{C}_R}\a{q_i}^2\right)^\frac{1}{2}\leq 2\en{\a{q_i}^2}^\frac{1}{2}\right. \\ \left.\textrm{and}\;\;\frac{1}{R}\left(\fint_{B_R}\a{\phi_i}^2+\a{\psi_i}^2+\a{\sigma_i}^2\right)^\frac{1}{2}+\frac{1}{R^2}\left(\fint_{\mathcal{C}_R}\a{\zeta_i}^2\right)^\frac{1}{2}\leq \frac{1}{C_0}\;\right\},\end{multline*}
which effectively defines the initial scale on which the $C^{1,\alpha}$-regularity of Proposition~\ref{i_excess_decay} begins to take effect, is expected to have stretched exponential moments in the sense of \cite[Theorem~1]{GNO4}.   Furthermore, again assuming a logarithmic Sobolev inequality, it should be possible to obtain a quantitative two-scale expansion for $a$-caloric functions like \cite[Corollary~3]{GNO4}.  Lastly, following the methods of \cite{FischerOtto}, it may be possible to prove higher order Liouville statements under a mild quantification of the ergodicity.

The paper is organized as follows.  The proofs are presented in the order of their appearance:   Theorem \ref{i_liouville}, Lemma~\ref{i_lemma1}, Proposition~\ref{i_sublinear}, Proposition~\ref{i_excess_decay} and Lemma~\ref{i_cacc}.  In order to simplify the notation, the statements and proofs are written for the non-symmetric scalar setting.  However, at the cost of increasing some constants, all of the arguments carry through unchanged for non-symmetric systems.  Throughout, the notation $\lesssim$ is used to denote a constant whose dependencies are specified in every case by the statement of the respective theorem, proposition or lemma.

\section{The proof of Theorem~\ref{i_liouville}}

Fix a coefficient field $a$ satisfying the conclusions of Lemma~\ref{i_lemma1}, Proposition~\ref{i_sublinear} and Proposition~\ref{i_excess_decay}, and suppose that $u$ is a distributional solution of
$$u_t=\nabla\cdot a\nabla u\;\;\textrm{on}\;\;\mathbb{R}^d\times(-\infty,0),$$
which is strictly subquadratic in the sense that, for some $\alpha\in(0,1)$,
$$\lim_{R\rightarrow\infty}\frac{1}{R^{1+\alpha}}\left(\fint_{\mathcal{C}_R}\a{u}^2\right)^{\frac{1}{2}}=0.$$
Fix $C_0=C_0(\alpha,d,\lambda)>0$ satisfying the hypothesis of Proposition~\ref{i_excess_decay}.  Then, using Proposition~\ref{i_sublinear}, fix $R_0>0$ such that, for every $R>R_0$, for each $i,j,k\in\{1,\ldots,d\}$,
$$\frac{1}{R}\left(\fint_{\mathcal{C}_R}\a{\phi_i-(\phi_i)_R}^2+\a{\psi_i-(\psi_i)_{t,R}}^2+\a{\sigma_{ijk}-(\sigma_{ijk})_{t,R}}^2\right)^{\frac{1}{2}}\leq\frac{1}{C_0},$$
and
$$\frac{1}{R^2}\left(\fint_{\mathcal{C}_R}\a{\zeta_{ij}-\zeta_{ij}(0)}^2\right)^\frac{1}{2}<\frac{1}{C_0},$$
and such that, for each $i\in\{1,\ldots,d\}$ and $R>R_0$,
$$\left(\fint_{\mathcal{C}_R}\a{q_i}^2\right)^\frac{1}{2}\leq 2\en{\a{q_i}^2}^\frac{1}{2}.$$

The definition of the excess, Proposition~\ref{i_excess_decay} and the Caccioppoli inequality (\ref{i_cacc_eq}) imply that, for each $R_0<\rho<R$,
$$\mathrm{Exc}(u;\rho)\lesssim \left(\frac{\rho}{R}\right)^{2\alpha}\mathrm{Exc}(u;R) \leq \left(\frac{\rho}{R}\right)^{2\alpha}\fint_{\mathcal{C}_R}\a{\nabla u}^2 \lesssim \frac{\rho^{2\alpha}}{(2R)^{2+2\alpha}}\fint_{\mathcal{C}_{2R}}\a{u}^2.$$
Therefore, since $u$ is strictly subquadratic,
$$\mathrm{Exc}(u;\rho)\lesssim \rho^{2\alpha}\limsup_{R\rightarrow\infty}\frac{1}{(2R)^{2+2\alpha}}\fint_{\mathcal{C}_{2R}}\a{u}^2=0.$$
This implies that, for every $\rho>0$,
$$\mathrm{Exc}(u;\rho)=0.$$
It is then immediate from the definition of parabolic excess (\ref{i_excess}), since $\mathcal{C}_{R_1}\subset\mathcal{C}_{R_2}$ whenever $R_1<R_2$, that there exists $\xi\in\mathbb{R}^d$ for which the difference
$$z(x,t):=u(x,t)-\xi\cdot x-\phi_\xi(x,t)$$
satisfies
$$\nabla z=0\;\;\textrm{in}\;\;\mathbb{R}^d\times(-\infty,0).$$
However, because $z$ is a distributional solution of
$$z_t=\nabla\cdot a\nabla z \;\;\textrm{in}\;\;\mathbb{R}^d\times(-\infty,0),$$
it follows that $z$ is necessarily constant in time as well.  Therefore, there exists $c\in\mathbb{R}$ such that $u=c+\xi\cdot x+\phi_\xi$, which completes the argument.

\section{The proof of Lemma~\ref{i_lemma1}}

The construction of the corrector $(\phi,\psi,\sigma,\zeta)$ will be achieved by lifting the relevant equations (\ref{i_corrector}), (\ref{i_psi}), (\ref{i_sigma}) and (\ref{i_zeta}) to the probability space $\Omega$, and thereby identifying $\phi$ by its stationary, finite energy gradient and time derivative, $\psi$ and $\sigma$ by their stationary, finite energy gradients, and $\zeta$ by its stationary time derivative.  For this, it is necessary to define the horizontal derivative of a random variable as induced by shifts of the coefficient field in space and time.  Then, these will be used to define an analogue of the Sobolev space $H^1$ on the probability space.

Given an $L^2(\Omega)$ random variable $f$, define, for each $i\in\{1,\ldots,d\}$, the horizontal derivative
\begin{equation}\label{phi_horizontal} D_if(a):=\lim_{h\rightarrow 0}\frac{f(a(\cdot+he_i,\cdot))-f(a)}{h},\end{equation}
where the above limit is understood in the sense of strong $L^2(\Omega)$-convergence.  Of course, it is not true in general that the above limit exists for every $f\in L^2(\Omega)$, but the horizontal derivatives are closed, densely defined operators on $L^2(\Omega)$, see \cite{papvar}, with domains, for $i\in\{1,\ldots,d\}$,
$$\mathcal{D}(D_i):=\{\;f\in L^2(\Omega)\;|\;D_if\;\;\textrm{exists as an element of}\;\;L^2(\Omega)\;\}.$$
Similarly, for each $f\in L^2(\Omega)$, define the horizontal time derivative
\begin{equation}\label{phi_time} D_0f(a):=\lim_{h\rightarrow 0}\frac{f(a(\cdot,\cdot+h))-f(a)}{h},\end{equation}
which is a closed, densely defined operator on $L^2(\Omega)$ with domain
$$\mathcal{D}(D_0):=\{\;f\in L^2(\Omega)\;|\;D_0f\;\;\textrm{exists as an element of}\;\;L^2(\Omega)\;\}.$$

The analogue of the Hilbert space $H^1$ is then defined on the probability space as the intersection
$$\mathcal{H}^1:=\cap_{i=0}^d\mathcal{D}(D_i),$$
equipped with the inner product, for each $f,g\in\mathcal{H}^1$,
$$(f,g)_{\mathcal{H}^1}=\en{fg}+\en{D_0f D_0g}+\en{Df\cdot Dg},$$
for the horizontal spatial gradient
$$Df:=(D_1f,\ldots, D_df).$$

Finally, define the space of spatial potentials
\begin{equation}\label{pot_1}L^2_{\textrm{pot}}(\Omega;\mathbb{R}^d):=\overline{\left\{\;Df\;|\;f\in\mathcal{H}^1\;\right\}}^{\;\;L^2(\Omega)-\textrm{weak}},\end{equation}
as the $L^2(\Omega)$-weak closure of spatial gradients arising from $\mathcal{H}^1$ functions.  Indeed, it is immediate from the weak convergence that elements of $L^2_{\textrm{pot}}(\Omega;\mathbb{R}^d)$ are potentials in the sense that every $F=(F_1,\ldots,F_d)\in L^2_{\textrm{pot}}(\Omega;\mathbb{R}^d)$ satisfies the distributional equality, for each $i,j\in\{1,\ldots,d\}$,
$$D_iF_j=D_jF_i.$$
In other words, every $F\in L^2_{\textrm{pot}}(\Omega;\mathbb{R}^d)$ is curl-free.

The following general fact about potential vector fields will be used in the construction of $\sigma$ and to prove the sublinearity for the corrector $(\phi, \psi, \sigma, \zeta)$.  It will be shown that, with respect to the sub-sigma-algebra of subsets that are invariant with respect to spatial translations of the coefficient fields, the conditional expectation of a potential vector field vanishes as a random variable.

\begin{lemma}\label{conditional_vanish} For every $F=(F_1,\ldots,F_d)\in L^2_\textrm{pot}(\Omega;\mathbb{R}^d)$, for each $i\in\{1,\ldots,d\}$,
\begin{equation}\label{cond_2} \en{F_i\;|\;\mathcal{F}_{\mathbb{R}^d}}=0\;\;\textrm{in}\;\;L^2(\Omega),\end{equation}
where $\en{\cdot\;|\;\mathcal{F}_{\mathbb{R}^d}}$ denotes the conditional expectation with respect to the sub-sigma-algebra $\mathcal{F}_{\mathbb{R}^d}\subset\mathcal{F}$ of subsets which are invariant with respect to spatial translations of the coefficient field.

In particular, for every $F=(F_1,\ldots,F_d)\in L^2_\textrm{pot}(\Omega;\mathbb{R}^d)$, for every $i\in\{1,\ldots,d\}$,
$$\en{F_i}=0.$$
\end{lemma}

The proof of Lemma~\ref{conditional_vanish} now follows.  Let $F=(F_1,\ldots,F_d)\in L^2_{\textrm{pot}}(\Omega;\mathbb{R}^d)$ be arbitrary.  Owing to the definition of the conditional expectation, it is sufficient to show that, for every $\mathcal{F}_{\mathbb{R}^d}$-measurable function $g\in L^2(\Omega)$, for each $i\in\{1,\ldots,d\}$,
\begin{equation}\label{cond_03}\en{F_ig}=0.\end{equation}
To prove (\ref{cond_03}), owing to definition (\ref{pot_1}) there exists a sequence of functions $\{\varphi_n\}_{n=1}^\infty\subset\mathcal{H}^1$ such that, as $n\rightarrow\infty$,
\begin{equation}\label{cond_1}D\varphi_n\rightharpoonup F\;\;\textrm{weakly in}\;\;L^2(\Omega;\mathbb{R}^d).\end{equation}
The weak convergence (\ref{cond_1}) implies that, for each $i\in\{1,\ldots,d\}$,
\begin{equation}\label{cond_3}\lim_{n\rightarrow\infty}\en{\left(D_i\varphi_n\right)g}=\en{F_i g}.\end{equation}
Then, for each $n\geq 1$, since $\varphi_n\in\mathcal{H}^1$, it follows that, since spatial translations of the coefficient field preserve the measure, see (\ref{i_stationary}), for each $i\in\{1,\ldots,d\}$,
\begin{equation}\label{cond_4}\en{\left(D_i\varphi_n\right)g}=\lim_{h\rightarrow 0}\frac{1}{h}\en{\left(\varphi_n(a(\cdot+he_i,\cdot))-\varphi_n(a)\right)g}=\lim_{h\rightarrow 0}\frac{1}{h}\en{\varphi_n\left(g(a(\cdot-he_i,\cdot))-g(a)\right)}=0,\end{equation}
where the final equality follows from the fact that $g\in L^2(\Omega)$ and the fact that $g$ is invariant with respect to spatial shifts of the coefficient field as an $\mathcal{F}_{\mathbb{R}^d}$-measurable function.  In combination, (\ref{cond_3}) and (\ref{cond_4}) imply (\ref{cond_03}).  Since the $\mathcal{F}_{\mathbb{R}^d}$-measurable $g\in L^2(\Omega)$ and $F\in L^2_{\textrm{pot}}(\Omega;\mathbb{R}^d)$ were arbitrary, this completes the proof of (\ref{cond_2}).

The final statement is then immediate since, for every $F\in L^2_{\textrm{pot}}(\Omega;\mathbb{R}^d)$ and $i\in\{1,\ldots,d\}$, the conditional expectation satisfies
$$\en{F_i}=\en{\en{\;F_i\;|\;\mathcal{F}_{\mathbb{R}^d}}}=0,$$
where the final equality follows from (\ref{cond_2}) and completes the proof of Lemma~\ref{conditional_vanish}.

\subsection{The construction of $\phi$}

For the construction of the corrector $\phi$, it is sufficient to construct, for each $k\in\{1,\ldots,d\}$, a stationary gradient $D\phi_k\in L^2_{\textrm{pot}}(\Omega;\mathbb{R}^d)$ and a stationary time-derivative $D_0\phi_k\in \mathcal{H}^{-1}$, where $\mathcal{H}^{-1}$ denotes the dual-space of $\mathcal{H}^1$, satisfying
\begin{equation}\label{phi_1}\en{D_0\phi_kf}+\en{Df\cdot a(D\phi_k+e_k)}=0\;\;\textrm{for every}\;\;f\in\mathcal{H}^1.\end{equation}
The corrector $\phi$ will then be defined on $\mathbb{R}^{d+1}$, for $\en{\cdot}$-a.e. $a$, by integration.

The first step is to introduce an approximation of (\ref{phi_1}) which is coercive with respect to the $\mathcal{H}^1$-norm.  The Riesz representation theorem and the uniform ellipticity of the ensemble (\ref{i_bounded}) guarantee that, for each $k\in\{1,\ldots,d\}$ and $\beta\in(0,1)$, there exists a unique element $\phi^\beta_k\in\mathcal{H}^1$ satisfying
\begin{equation}\label{phi_2}\beta\en{\phi^\beta_kf}+\beta\en{D_0\phi^\beta_kD_0f}+\en{D_0\phi^\beta_kf}+\en{Df\cdot a(D\phi^\beta_k+e_k)}=0\;\;\textrm{for every}\;\;f\in\mathcal{H}^1.\end{equation}
Therefore, for each $k\in\{1,\ldots,d\}$ and $\beta\in(0,1)$, after testing (\ref{phi_2}) with $\phi^\beta_k$, the uniform ellipticity of the ensemble and H\"older's inequality yield the estimate
\begin{equation}\label{phi_3}\en{\a{D\phi^\beta_k}^2}+\beta\en{\a{\phi^\beta_k}^2+\a{D_0\phi^\beta_k}^2}\lesssim 1,\end{equation}
where the fact that
\begin{equation}\label{phi_303}\en{D_0\phi^\beta_k\phi^\beta_k}=\frac{1}{2}\en{D_0(\phi^\beta_k)^2}=\lim_{h\rightarrow 0}\frac{1}{2h}\en{(\phi^\beta_k)^2(a(\cdot,\cdot+h))-(\phi^\beta_k)^2(a)}=0\end{equation}
is also used to obtain (\ref{phi_3}), and follows from the fact that shifts of the coefficient field in time and space preserve the underlying measure of the ensemble, see (\ref{i_stationary}), and $\phi^\beta_k\in L^2(\Omega)$.

Then, for each $k\in\{1,\ldots,d\}$ and $\beta\in(0,1)$, equation (\ref{phi_2}), H\"older's inequality and (\ref{phi_3}) imply that, for each $f\in\mathcal{H}^1$,
\begin{equation}\label{phi_4}\a{\en{D_0\phi^\beta_kf}}\lesssim \norm{f}_{\mathcal{H}^1}\;\;\textrm{and, therefore,}\;\;\norm{D_0\phi^\beta_k}_{\mathcal{H}^{-1}}\lesssim 1.\end{equation}
Hence, for each $k\in\{1,\ldots,d\}$, the definition of the potential space (\ref{pot_1}) with estimates (\ref{phi_3}) and (\ref{phi_4}) imply that there exist $\Psi_k\in L^2_{\textrm{pot}}(\Omega;\mathbb{R}^d)$ and $\xi_k\in \mathcal{H}^{-1}$ such that, after passing to a subsequence $\{\beta_j\rightarrow 0\}_{j=1}^\infty$, as $j\rightarrow\infty$,
\begin{equation}\label{phi_5}D\phi^{\beta_j}_k\rightharpoonup \Psi_k\;\;\textrm{weakly in}\;\;L^2(\Omega;\mathbb{R}^d)\;\;\textrm{and}\;\;D_0\phi^{\beta_j}_k\rightharpoonup \xi_k\;\;\textrm{weakly in}\;\;\mathcal{H}^{-1}.\end{equation}
The convergence (\ref{phi_5}) combined with equation (\ref{phi_2}) and estimate (\ref{phi_3}) prove that, for each $k\in\{1,\ldots,d\}$,
$$\en{\xi_k f}+\en{Df\cdot a(\Psi_k+e_k)}=0\;\;\textrm{for every}\;\;f\in\mathcal{H}^1.$$
Finally, for each $k\in\{1,\ldots,d\}$, as weak limits of functions $\phi^\beta_k\in \mathcal{H}^1$, the pair $(\Psi_k,\xi_k)$ are curl-free in the sense that, for each $i\in\{1,\ldots,d\}$, distributionally
$$D_0\Psi_{ki}=D_i\xi_k.$$
Therefore, for each $k\in\{1,\ldots,d\}$, the argument is completed by defining $D\phi_k:=\Psi_k$ and $D_0\phi_k:=\xi_k$.

\subsection{The construction of $\psi$}

For each $k\in\{1,\ldots,d\}$, let $D\phi_k\in L^2_{\textrm{pot}}(\Omega;\mathbb{R}^d)$ denote the stationary gradient corresponding to $\phi_k$, which was constructed in the previous step.  Then, for each $k\in\{1,\ldots,d\}$, define the lift of the flux $q_k$ to the probability space according to the rule
\begin{equation}\label{psi_flux}Q_k:=a(D\phi_k+e_k).\end{equation}
The existence of $\psi$ follows from the following general fact.

\begin{lemma}\label{psi_construction}  For every $F\in L^2(\Omega;\mathbb{R}^d)$, there exists $\Psi \in L^2_{\textrm{pot}}(\Omega;\mathbb{R}^d)$ satisfying
\begin{equation}\label{psi_1} \en{\Psi\cdot Df}=\en{F\cdot Df}\;\;\textrm{for every}\;\;f\in\mathcal{H}^1.\end{equation}
\end{lemma}

The existence of $\psi$ follows from Lemma~\ref{psi_construction} in the following way.  For each $k\in\{1,\ldots,d\}$, choose $F=Q_k$ and define $D\psi_k:=\Psi$, which defines $\psi_k$, for $\en{\cdot}$-a.e. $a$, as a function on $\mathbb{R}^d$ via integration.  Then, for $\en{\cdot}$-a.e. $a$, for each $k\in\{1,\ldots,d\}$, the function $\psi_k$ is extended to $\mathbb{R}^{d+1}$ as a stationary function in time.

In order to prove Lemma~\ref{psi_construction}, the Riesz representation theorem asserts that, for each $\beta\in(0,1)$, there exists a unique $\psi^\beta\in\mathcal{H}^1$ satisfying
\begin{equation}\label{psi_2} \beta\en{\psi^\beta f}+\beta\en{D_0\psi^\beta D_0f}+\en{D\psi^\beta\cdot Df}=\en{F\cdot Df}\;\;\textrm{for every}\;\;f\in\mathcal{H}^1.\end{equation}
For each $\beta\in(0,1)$, after testing (\ref{psi_2}) with $\psi^\beta$, H\"older's inequality and Young's inequality yield the estimate
\begin{equation}\label{psi_3}\en{\a{D\psi^\beta}^2}+\beta\en{\a{\psi^\beta}^2+\a{D_0\psi^\beta}^2}\lesssim 1.\end{equation}
Therefore, the definition of the potential space (\ref{pot_1}) and estimate (\ref{psi_3}) imply that there exists $\Psi\in L^2_{\textrm{pot}}(\Omega;\mathbb{R}^d)$ such that, after passing to a subsequence $\{\beta_j\rightarrow 0\}_{j=1}^\infty$, as $j\rightarrow\infty$,
\begin{equation}\label{psi_4}D\psi^{\beta_j}\rightharpoonup \Psi\;\;\textrm{weakly in}\;\;L^2(\Omega;\mathbb{R}^d).\end{equation}
In combination with equation (\ref{psi_2}), estimates (\ref{psi_3}) and the convergence (\ref{psi_4}) imply that $\Psi$ satisfies (\ref{psi_1}), which completes the proof of Lemma~\ref{psi_construction}.

\subsection{The construction of $\sigma$}

For each $k\in\{1,\ldots,d\}$, let $D\psi_k\in L^2_{\textrm{pot}}(\Omega;\mathbb{R}^d)$ denote the stationary gradient corresponding to $\psi_k$ constructed in the previous step and let $Q_k$ denote the lift of the flux from (\ref{psi_flux}).  Lemma~\ref{psi_construction} applies directly to this situation, and proves that, for each $i,j,k\in\{1,\ldots,d\}$, there exists $\Sigma_{ijk}\in L^2_{\textrm{pot}}(\Omega;\mathbb{R}^d)$ satisfying
\begin{equation}\begin{aligned}\label{sigma_1} \en{\Sigma_{ijk}\cdot Df}= & \en{(Q_i-D\psi_i)_k D_jf}-\en{(Q_i-D\psi_i)_jD_kf} \\ = & \en{\left((Q_i-D\psi_i)_ke_j-(Q_i-D\psi_i)_je_k\right)\cdot Df}\;\;\textrm{for every}\;\;f\in\mathcal{H}^1.\end{aligned}\end{equation}
Then, for each $i,j,k\in\{1,\ldots,d\}$, the definition $D\sigma_{ijk}:=\Sigma_{ijk}$ defines $\sigma_{ijk}$, for $\en{\cdot}$-a.e. $a$, on $\mathbb{R}^d$ via integration.  These functions are then extended to $\mathbb{R}^{d+1}$, for each $i,j,k\in\{1,\ldots,d\}$, for $\en{\cdot}$-a.e. $a$, as stationary functions in time.

Since it is clear from the proof of existence that, for each $i,j,k\in\{1,\ldots,d\}$, the gradients can be constructed to satisfy
$$\Sigma_{ijk}=-\Sigma_{ikj},$$
after integrating it follows that, for each $i,j,k\in\{1,\ldots,d\}$,
$$\sigma_{ijk}=-\sigma_{ikj}.$$
Or, perhaps more simply, for each $i\in\{1,\ldots,d\}$, one may first construct $\sigma_{ijk}$ for every $j>k\in\{1,\ldots,d\}$ and then simply define $\sigma_{ikj}:=-\sigma_{ijk}$.

It remains to prove that, for each $i\in\{1,\ldots,d\}$,
\begin{equation}\label{sigma_2} D\cdot \sigma_i=(Q_i-D\psi_i)-\en{Q_i\;|\;\mathcal{F}_{\mathbb{R}^d}},\end{equation}
where, for each $i,j\in\{1,\ldots,d\}$,
\begin{equation}\label{sigma_20}(D\cdot\sigma_i)_j=\sum_{k=1}^dD_k\sigma_{ijk}.\end{equation}

To simplify the notation in what follows, define the vector, for each $i\in\{1,\ldots,d\}$,
\begin{equation}\label{sigma_21}\Psi_i:=(Q_i-D\psi_i),\end{equation}
where the construction of $\psi_i$ guarantees that the vector $\Psi_i$ is divergence-free in the sense of the distributional equality
\begin{equation}\label{sigma_3}D_l\Psi_{il}=0.\end{equation}
Equation (\ref{sigma_2}) now follows from the following distributional equalities.  For each $i,j\in\{1,\ldots,d\}$, thanks to equation (\ref{sigma_1}), in the sense of distributions
$$D_l(D_l(D\cdot\sigma_i)_j)=D_l(D_l(D_k\sigma_{ijk}) =D_k(D_l(D_l\sigma_{ijk}))=D_k(D_k\Psi_{ij}-D_j\Psi_{ik}).$$
Therefore, for each $i,j\in\{1,\ldots,d\}$, in view of (\ref{sigma_3}) and after relabeling the final integral,
$$D_l(D_l(D\cdot\sigma_i)_j)=D_k(D_k\Psi_{ij}-D_j\Psi_{ik})=D_k(D_k\Psi_{ij})-D_j(D_k\Psi_{ik})=D_l(D_l\Psi_{ij}).$$
Hence, in the sense of distributions, for each $i,j\in\{1,\ldots,d\}$,
\begin{equation}\label{sigma_4}D_l\left(D_l\left((D\cdot\sigma_i)_j-\Psi_{ij}\right)\right)=0.\end{equation}

Equation (\ref{sigma_4}) implies that, for each $i,j\in\{1,\ldots,d\}$, the difference $\left((D\cdot\sigma_i)_j-\Psi_{ij}\right)$ is invariant with respect to spatial translations of the coefficient field.  That is, for each $i,j\in\{1,\ldots,d\}$,
\begin{equation}\label{sigma_5}\left((D\cdot\sigma_i)_j-\Psi_{ij}\right)=\en{\left((D\cdot\sigma_i)_j-\Psi_{ij}\right)\;|\;\mathcal{F}_{\mathbb{R}^d}}.\end{equation}
The fact that (\ref{sigma_5}) implies (\ref{sigma_2}) follows from Lemma~\ref{conditional_vanish}.  Indeed, for each $i\in\{1,\ldots,d\}$,
$$D\psi_i \in L^2_{\textrm{pot}}(\Omega;\mathbb{R}^d),$$
and by a straightforward repetition of the arguments leading to Lemma~\ref{conditional_vanish}, for each $i\in\{1,\ldots,d\}$,
$$\en{\;(D\cdot\sigma_i)\;|\;\mathcal{F}_{\mathbb{R}^d}\;}=0\;\;\textrm{in}\;\; L^2(\Omega;\mathbb{R}^d).$$
Therefore, for each $i,j\in\{1,\ldots,d\}$,
$$\left((D\cdot\sigma_i)_j-\Psi_{ij}\right)=\en{\left((D\cdot\sigma_i)_j-\Psi_{ij}\right)\;|\;\mathcal{F}_{\mathbb{R}^d}}=-\en{Q_{ij}\;|\;\mathcal{F}_{\mathbb{R}^d}},$$
which is (\ref{sigma_2}).  This completes the argument.

\subsection{The construction of $\zeta$}

The construction of $\zeta$ is explicit.  Namely, for each $i\in\{1,\ldots,d\}$, for the lift of the flux $Q_i$ defined in (\ref{psi_flux}), define the stationary derivative in time according to the rule
$$D_0\zeta_i=\en{Q_i\;|\;\mathcal{F}_{\mathbb{R}^d}}-\en{Q_i}.$$
For each $i\in\{1,\ldots,d\}$, the function $\zeta_i$ is defined on $\mathbb{R}$, for $\en{\cdot}$-a.e $a$, by integration in time and extended to $\mathbb{R}^{d+1}$ as a spatially constant function.

\subsection{The boundedness and uniform ellipticity of $\ah$}

For the reader's convenience, the argument of \cite[Lemma~2]{GNO4} is repeated here.  For each $\xi\in\mathbb{R}^d$, the linearity and (\ref{i_homogenized_coefficients}) assert that the homogenized coefficients are defined according to the rule
$$\ah\xi:=\en{a(\nabla\phi_\xi+\xi)}.$$
It is first shown that, for each $\xi\in\mathbb{R}^d$,
\begin{equation}\label{ah_1} \a{\ah \xi}\leq\frac{1}{\lambda} \a{\xi}.\end{equation}
For each $\xi\in\mathbb{R}^d$, since $\phi_\xi$ satisfies (\ref{i_corrector}), the uniform ellipticity of the ensemble (\ref{i_bounded}) and Jensen's inequality imply that
\begin{multline*}\a{\ah\xi}^2=\a{\en{a(\nabla\phi_\xi+\xi)}}^2\leq \\ \en{\a{a(\nabla\phi_\xi+\xi)}^2}\leq\en{\a{(\nabla\phi_\xi+\xi)}^2}\leq \frac{1}{\lambda}\en{(\nabla\phi_\xi+\xi)\cdot a(\nabla\phi_\xi+\xi)}.\end{multline*}
Then, using the corrector equation (\ref{i_corrector}), (\ref{phi_303}) and the Cauchy-Schwarz inequality,
$$\a{\ah\xi}^2\leq \frac{1}{\lambda}\en{(\nabla\phi_\xi+\xi)\cdot a(\nabla\phi_\xi+\xi)}=\frac{1}{\lambda}\xi\cdot\en{a(\nabla\phi_\xi+\xi)}\leq \frac{1}{\lambda}\a{\xi}\a{\ah\xi}.$$
Dividing by $\a{\ah\xi}$ yields (\ref{ah_1}) and completes the proof.

It remains only to prove that, for each $\xi\in\mathbb{R}^d$,
\begin{equation}\label{ah_2}\lambda\a{\xi}^2\leq \xi\cdot \ah\xi.\end{equation}
This follows from the convexity of the map
\begin{equation}\label{convex_map}(X,v)\in\mathcal{S}(d)_{> 0}\times\mathbb{R}^d\rightarrow v\cdot X^{-1}v,\end{equation}
where $\mathcal{S}(d)_{>0}$ denotes the space of positive, $d\times d$ symmetric matrices.  Indeed, if $X\in\mathcal{S}(d)_{>0}$ and $v\in\mathbb{R}^d$,
$$\frac{1}{2}\left(v\cdot X^{-1}v\right)=\sup_{w\in\mathbb{R}^d}\left\{w\cdot v-\frac{1}{2}w\cdot Xw\right\},$$
is a supremum over linear functions in $(X,v)$, and therefore the map (\ref{convex_map}) is convex.  Hence, for each $\xi\in\mathbb{R}^d$, using the corrector equation (\ref{i_corrector}) and Jensen's inequality,
$$\begin{aligned} \xi\cdot \ah \xi = & \en{(\nabla\varphi_\xi+\xi)\cdot a(\nabla\varphi_\xi+\xi)}=\en{(\nabla\phi_\xi+\xi)\cdot a_{\textrm{sym}}(\nabla\phi_\xi+\xi)} \\ = & \en{(\nabla\phi_\xi+\xi)\cdot (a_{\textrm{sym}}^{-1})^{-1}(\nabla\phi_\xi+\xi)}\geq \en{\nabla\phi_\xi+\xi}\cdot \en{(a_{\textrm{sym}}^{-1})}^{-1}\en{\nabla\phi_\xi+\xi} \\ \geq & \lambda\a{\xi}^2,\end{aligned}$$
where $a_{\textrm{sym}}$ denotes the symmetric part of $a$, and where the final inequality is obtained using the boundedness of the ensemble (\ref{i_bounded}) and the vanishing expectation of the gradient from Lemma~\ref{conditional_vanish}.  This completes the proof of (\ref{ah_2}).

\section{The proof of Proposition~\ref{i_sublinear}}

\subsection{The sublinearity of $\sigma$ and $\psi$}

The sublinearity of $\sigma$ and $\psi$ will follow from the following general fact.  Recall that, for a function $\varphi:\mathbb{R}^{d+1}\rightarrow\mathbb{R}$, for each $R>0$ and $t\in\mathbb{R}$,
$$(\varphi)_R=\fint_{\mathcal{C}_R}\varphi\;\;\textrm{and}\;\;(\varphi)_{t,R}=\fint_{B_R}\varphi(\cdot,t).$$

\begin{lemma}\label{lemma_sigma_sub}Suppose that $\varphi$ is a scalar random field on $\mathbb{R}^{d+1}$ which is stationary in time and has a stationary, finite energy gradient $\nabla\varphi$ in the potential space $L^2_{\textrm{pot}}(\Omega;\mathbb{R}^d)$.  That is, assume that
\begin{equation}\label{new_1000}\en{\a{\nabla\varphi}^2}<\infty\;\;\textrm{with}\;\;\nabla\varphi\in L^2_{\textrm{pot}}(\Omega;\mathbb{R}^d).\end{equation}
Then, for $\en{\cdot}$-a.e. $a$, the normalized large-scale $L^2$-averages of $\varphi$ are strictly sublinear in the sense that
\begin{equation}\label{new_1}\lim_{R\rightarrow\infty}\frac{1}{R^2}\fint_{\mathcal{C}_R}\a{\varphi-(\varphi)_{t,R}}^2=0.\end{equation}
\end{lemma}

To prove Lemma~\ref{lemma_sigma_sub}, since $\varphi$ is stationary in time, it will first be shown that, for $\en{\cdot}$-a.e. $a$, the normalized large-scale $L^2$-averages of $\varphi(\cdot,0)$ on large balls are sublinear.  Define, for each $\epsilon>0$,
$$\varphi^\epsilon(\cdot):=\epsilon\varphi(\frac{\cdot}{\epsilon},0).$$
Then, for $\en{\cdot}$-a.e. $a$,
\begin{equation}\label{new_3}\lim_{R\rightarrow\infty}\frac{1}{R^2}\fint_{B_R}\a{\varphi(\cdot,0)-(\varphi)_{0,R}}^2=\lim_{\epsilon\rightarrow 0} \fint_{B_1}\a{\varphi^\epsilon-(\varphi^\epsilon)}^2=0,\end{equation}
where the first equality is immediate by scaling if, for each $\epsilon>0$,
$$(\varphi^\epsilon):=\fint_{B_1}\varphi^\epsilon.$$

To obtain (\ref{new_3}), the Poincar\'e inequality in space, together with the ergodic theorem, imply that, for $\en{\cdot}$-a.e. $a$,
$$\limsup_{\epsilon\rightarrow 0}\fint_{B_1}\a{\varphi^\epsilon-(\varphi^\epsilon)}^2\leq \limsup_{\epsilon\rightarrow 0}\fint_{B_1}\a{\nabla\varphi^\epsilon}^2=\en{\a{\nabla\varphi}^2\;|\;\mathcal{F}_{\mathbb{R}^d}}<\infty,$$
where the final term is finite, for $\en{\cdot}$-a.e. $a$, thanks to (\ref{new_1000}).  Therefore, for $\en{\cdot}$-a.e. $a$, the sequence $\{(\varphi^\epsilon-(\varphi^\epsilon))\}_{\epsilon\in(0,1)}$ is bounded in $H^1(B_1)$, and hence compact in $L^2(B_1)$ by the Rellich-Kondrachov embedding theorem.  Since the ergodic theorem, see Krengel \cite[Theorem~2.3]{Krengel}, implies that, for $\en{\cdot}$-a.e. $a$, as $\epsilon\rightarrow 0$, the gradient
$$\nabla\varphi^\epsilon\rightharpoonup \en{\nabla\varphi\;|\;\mathcal{F}_{\mathbb{R}^d}}=0\;\;\textrm{weakly in}\;\;L^2(B_1;\mathbb{R}^d),$$
where the vanishing of the conditional expectation follows from Lemma~\ref{conditional_vanish} and (\ref{new_1000}), the compact embedding implies that, for $\en{\cdot}$-a.e $a$,
$$(\varphi^\epsilon-(\varphi^\epsilon))\rightarrow 0\;\;\textrm{strongly in}\;\;L^2(B_1),$$
which proves (\ref{new_3}).

The sublinearity of the normalized large-scale averages at time zero will now be upgraded to the sublinearity of the normalized large-scale averages on parabolic cylinders using Egorov's theorem.  Precisely, using (\ref{new_3}), it will be shown that, for $\en{\cdot}$-a.e. $a$,
\begin{equation}\label{new_8}\lim_{R\rightarrow\infty}\frac{1}{R}\left(\fint_{\mathcal{C}_R}\a{\varphi-(\varphi)_{t,R}}^2\right)^{\frac{1}{2}}=0.\end{equation}

To prove (\ref{new_8}), for $\epsilon\in(0,1)$ use Egorov's theorem and (\ref{new_3}) to find a measurable subset $A_\epsilon\subset\Omega$ and $R_\epsilon>0$ such that, for every $a\in A_\epsilon$ and $R>R_\epsilon$,
\begin{equation}\label{new_2}\frac{1}{R}\left(\fint_{B_R}\a{\varphi(\cdot,0)-(\varphi)_{0,R}}^2\right)^{\frac{1}{2}}<\epsilon\;\;\textrm{with}\;\;\en{\chi_{A_\epsilon}}>1-\epsilon,\end{equation}
where $\chi_{A_\epsilon}\in L^\infty(\Omega)$ denotes the indicator function of $A_\epsilon$.  The large-scale averages on cylinders appearing in (\ref{new_8}) will be decomposed according to the event $A_\epsilon$.

For each $R>0$ and $\epsilon\in(0,1)$, use the triangle inequality to form the decomposition
\begin{multline}\label{new_6}\frac{1}{R}\left(\fint_{\mathcal{C}_R}\a{\varphi-(\varphi)_{t,R}}^2\right)^{\frac{1}{2}}\leq \frac{1}{R}\left(\fint_{-R^2}^0\chi_{A_\epsilon}(a(\cdot,\cdot+t))\fint_{B_R}\a{\left(\varphi-(\varphi)_{t,R}\right)}^2\right)^{\frac{1}{2}} \\ + \frac{1}{R}\left(\fint_{-R^2}^0\left(1-\chi_{A_\epsilon}(a(\cdot,\cdot+t))\right)\fint_{B_R}\a{\left(\varphi-(\varphi)_{t,R}\right)}^2\right)^{\frac{1}{2}}.\end{multline}
For the second term of (\ref{new_6}), for each $R>0$ and $\epsilon\in(0,1)$, H\"older's inequality and the Sobolev embedding theorem imply that
\begin{multline}\label{new_6000}\frac{1}{R}\left(\fint_{-R^2}^0\left(1-\chi_{A_\epsilon}(a(\cdot,\cdot+t))\right)\fint_{B_R}\a{\left(\varphi-(\varphi)_{t,R}\right)}^2\right)^{\frac{1}{2}} \lesssim\\ \frac{1}{R}\left(\fint_{\mathcal{C}_R}\left(1-\chi_{A_\epsilon}(a(\cdot,\cdot+t))\right)\right)^\frac{1}{d}\left(\fint_{\mathcal{C}_R}\a{\left(\varphi-(\varphi)_{t,R}\right)}^\frac{2d}{d-2}\right)^\frac{d-2}{2d}\lesssim \\ \left(\fint_{\mathcal{C}_R}\left(1-\chi_{A_\epsilon}(a(\cdot,\cdot+t))\right)\right)^\frac{1}{d}\left(\fint_{\mathcal{C}_R}\a{\nabla\varphi}^2\right)^\frac{1}{2},\end{multline}
where the argument is written only for the case $d\geq 3$, since the modifications necessary to handle the cases $d=1$ and $d=2$ are straightforward and follow from the Sobolev embedding theorem.

After combining (\ref{new_6}) and (\ref{new_6000}), it is then immediate from the ergodic theorem and the definition of $A_\epsilon$ from (\ref{new_2}) that, for each $\epsilon>0$, for $\en{\cdot}$-a.e. $a$,
\begin{equation}\label{new_5}\limsup_{R\rightarrow\infty}\frac{1}{R}\left(\fint_{\mathcal{C}_R}\a{\varphi-(\varphi)_{t,R}}^2\right)^{\frac{1}{2}}\leq \en{(1-\chi_{A_\epsilon})\;|\;\mathcal{F}_t}^\frac{1}{d}\en{\a{\nabla\varphi}^2}^\frac{1}{2}+\epsilon,\end{equation}
where $\en{\cdot\;|\;\mathcal{F}_t}$ denotes the conditional expectation with respect to the sub-sigma-algebra $\mathcal{F}_t\subset\mathcal{F}$ of subsets of $\Omega$ which are invariant with respect to translations of the coefficient fields in time.

Since, as $\epsilon\rightarrow 0$, the choice (\ref{new_2}) implies that
$$\en{(1-\chi_{A_\epsilon})\;|\;\mathcal{F}_t} \rightarrow0\;\;\textrm{strongly in}\;\;L^1(\Omega),$$
using (\ref{new_1000}), after choosing a countable sequence $\{\epsilon_n\rightarrow 0\}_{n=1}^\infty$, it follows from (\ref{new_6}) and (\ref{new_5}) that, for $\en{\cdot}$-a.e. $a$,
\begin{equation}\label{new_7}\lim_{R\rightarrow\infty}\frac{1}{R}\left(\fint_{\mathcal{C}_R}\a{\varphi-(\varphi)_{t,R}}^2\right)^{\frac{1}{2}}=0,\end{equation}
which proves (\ref{new_8}).  This completes the proof of Lemma~\ref{lemma_sigma_sub}, and thereby proves the sublinearity of $\sigma$ and $\psi$ thanks to Lemma~\ref{i_lemma1}.

\subsection{The sublinearity of $\phi$}

The sublinearity of $\phi$ will follow from the following general fact.

\begin{lemma}\label{lemma_sub_phi}Suppose that $\varphi$ is a scalar random field on $\mathbb{R}^{d+1}$ which has a stationary, finite energy gradient $\nabla\varphi$ in the potential space $L^2_{\textrm{pot}}(\Omega;\mathbb{R}^d)$.  That is, assume that
\begin{equation}\label{new_0}\en{\a{\nabla\varphi}^2}<\infty\;\;\textrm{with}\;\;\nabla\varphi\in L^2_{\textrm{pot}}(\Omega;\mathbb{R}^d).\end{equation}
Furthermore, assume that, for $\en{\cdot}$-a.e. $a$, the field $\varphi$ satisfies
\begin{equation}\label{sub_0}\varphi_t=\nabla\cdot F\;\;\textrm{distributionally in}\;\;\mathbb{R}^{d+1},\end{equation}
for the stationary extension of a finite energy field $F\in L^2(\Omega;\mathbb{R}^d)$.  Then, for $\en{\cdot}$-a.e. $a$, the normalized large-scale $L^2$-averages of $\varphi$ are strictly sublinear in the sense that
\begin{equation}\label{sub_5} \lim_{R\rightarrow \infty}\frac{1}{R}\left(\fint_{\mathcal{C}_R}\a{\varphi-(\varphi)_R}^2\right)^\frac{1}{2}=0.\end{equation}
\end{lemma}

To prove Lemma~\ref{lemma_sub_phi}, for each $R>0$, use the triangle inequality to obtain
\begin{equation}\label{sub_7}\frac{1}{R}\left(\fint_{\mathcal{C}_R}\a{\varphi-(\varphi)_R}^2\right)^{\frac{1}{2}}\leq \frac{1}{R}\left(\fint_{\mathcal{C}_R}\a{\varphi-(\varphi)_{t,R}}^2\right)^\frac{1}{2}+\frac{1}{R}\left(\fint_{\mathcal{C}_R}\a{(\varphi)_{t,R}-(\varphi)_R}^2\right)^\frac{1}{2}.\end{equation}
Since $\varphi$ has a stationary spatial gradient, for each $R>0$ and $t\in\mathbb{R}$, the scalar random field
$$(x,t)\rightarrow(\varphi-(\varphi)_{t,R})(x,t)=\fint_{B_R}\varphi(x,t)-\varphi(y,t)\dy$$
is stationary in time, with a stationary, finite-energy gradient in the potential space.  Therefore, Lemma~\ref{lemma_sigma_sub} applies to this random field, and asserts that, for $\en{\cdot}$-a.e. $a$,
\begin{equation}\label{sub_8}\lim_{R\rightarrow\infty}\frac{1}{R}\left(\fint_{\mathcal{C}_R}\a{\varphi-(\varphi)_{t,R}}^2\right)^\frac{1}{2}=0.\end{equation}

It remains to prove that, for $\en{\cdot}$-a.e. $a$,
$$\lim_{R\rightarrow 0}\frac{1}{R}\left(\fint_{\mathcal{C}_R}\a{(\varphi)_{t,R}-(\varphi)_R}^2\right)^\frac{1}{2}=\lim_{R\rightarrow 0}\frac{1}{R}\left(\fint^0_{-R^2}\a{(\varphi)_{t,R}-(\varphi)_R}^2\right)^\frac{1}{2}=0.$$
Let $\rho\in\C^\infty_c(\mathbb{R}^d)$ be a smooth, symmetric convolution kernel satisfying $\supp(\rho)\subset B_1$ and, for each $\epsilon>0$, define the rescaling $\rho^\epsilon(\cdot)=\epsilon^{-d}\rho(\frac{\cdot}{\epsilon})$.  Then, for each $\epsilon>0$, define the spatial convolution, for each $x\in\mathbb{R}^d$ and $t\in\mathbb{R}$,
$$\varphi^\epsilon(x,t):=\int_{\mathbb{R}^d} \rho^\epsilon(y-x)\varphi(y,t)\dy.$$
The introduction of the convolution kernel provides a test function which will be used to apply the equation (\ref{sub_0}) satisfied by $\varphi$.

First, for each $R>0$ and $\epsilon\in(0,R)$, it follows from the support of the convolution kernel, Fubini's theorem and Jensen's inequality that
\begin{equation}\label{sub_9}\begin{aligned}\a{(\varphi)_R-(\varphi^\epsilon)_R}^2 = & \a{\fint_{\mathcal{C}_R}\int_{\mathbb{R}^d} \rho^\epsilon(x)\left(\varphi(x+y,t)-\varphi(y,t)\right)\dx\dy\dt}^2 \\ \leq & \epsilon^2\a{\int_0^1\int_{\mathbb{R}^d}\rho^\epsilon(x)\fint_{\mathcal{C}_R}\a{\nabla\varphi(y+rx,t)}\dy\dx\dt\dr}^2 \\ \leq & \epsilon^2\int_0^1\int_{\mathbb{R}^d}\rho^\epsilon(x)\fint_{\mathcal{C}_R}\a{\nabla\varphi(y+rx,t)}^2\dy\dx\dt\dr \\  \lesssim & \epsilon^2\fint_{\mathcal{C}_{2R}}\a{\nabla\varphi}^2.\end{aligned}\end{equation}
Similarly, the identical argument yields, for each $R>0$, $t\in\mathbb{R}$ and $\epsilon\in(0,R)$,
\begin{equation}\label{sub_900}\a{(\varphi)_{t,R}-(\varphi^\epsilon)_{t,R}}^2\lesssim \epsilon^2\fint_{B_{2R}}\a{\nabla\varphi(\cdot,t)}^2.\end{equation}

Therefore, for each $R>0$ and $\epsilon\in(0,R)$, the triangle inequality, (\ref{sub_9}) and (\ref{sub_900}) imply that, after adding and subtracting $((\varphi^\epsilon)_{t,R}-(\varphi^\epsilon)_R)$,
\begin{equation}\label{sub_10}\frac{1}{R}\left(\fint^0_{-R^2}\a{(\varphi)_{t,R}-(\varphi)_R}^2\right)^{\frac{1}{2}} \lesssim \frac{1}{R}\left(\fint^0_{-R^2}\a{(\varphi^\epsilon)_{t,R}-(\varphi^\epsilon)_R}^2\right)^{\frac{1}{2}}+\frac{\epsilon}{R}\left(\fint_{\mathcal{C}_{2R}}\a{\nabla\varphi}^2\right)^\frac{1}{2}.\end{equation}
For the first term on the right hand side of (\ref{sub_10}), for each $R\geq 1$ and $\epsilon\in(0,R)$, the equation (\ref{sub_0}) satisfied by $\varphi$ and the Poincar\'e inequality in time yields
\begin{equation}\label{sub_11}\begin{aligned} \frac{1}{R}\left(\fint^0_{-R^2}\a{(\varphi^\epsilon)_{t,R}-(\varphi^\epsilon)_R}^2\right)^\frac{1}{2} \lesssim & \left(\fint^0_{-R^2}\a{\partial_t(\varphi^\epsilon)_{t,R}}^2\right)^\frac{1}{2} \\ \lesssim & \left(\fint^0_{-R^2}\a{\fint_{B_R}\int_{\mathbb{R}^d}\nabla\rho^\epsilon(x-y)\cdot F(x,t)\dx\dy}^2\dt\right)^\frac{1}{2} \\ \lesssim & \left(\fint^0_{-R^2}\a{\int_{\mathbb{R}^d}\a{\nabla\rho^\epsilon(x)}\fint_{B_R(x)}\a{F(y,t)}\dy\dx}^2\dt\right)^\frac{1}{2} \\ \lesssim & \left(\fint^0_{-R^2}\a{\left(\int_{\mathbb{R}^d}\a{\nabla\rho^\epsilon(x)}\;\dx\right)\left(\fint_{B_{2R}}\a{F(y,t)}\;\dy\right)}^2\dt\right)^\frac{1}{2}.\end{aligned}\end{equation}
Then, continuing with (\ref{sub_11}), the definition of the convolution kernel and Jensen's inequality yield
\begin{equation}\label{sub_12}\begin{aligned} \frac{1}{R}\left(\fint^0_{-R^2}\a{(\varphi^\epsilon)_{t,R}-(\varphi^\epsilon)_R}^2\right)^\frac{1}{2} \lesssim &  \frac{1}{\epsilon}\left(\fint^0_{-R^2}\left(\fint_{B_{2R}}\a{F(y,t)}\;\dy\right)^2\dt\right)^\frac{1}{2} \\ \lesssim & \frac{1}{\epsilon}\left(\fint_{\mathcal{C}_{2R}}\a{F}^2\right)^\frac{1}{2}. \end{aligned}\end{equation}

Therefore, combining (\ref{sub_10}) with (\ref{sub_12}), for each $R>0$ and $\epsilon\in(0,R)$,
$$\frac{1}{R}\left(\fint^0_{-R^2}\a{(\varphi)_{t,R}-(\varphi)_R}^2\right)^{\frac{1}{2}}\lesssim \frac{1}{\epsilon}\left(\fint_{\mathcal{C}_{2R}}\a{F}^2\right)^\frac{1}{2}+\frac{\epsilon}{R}\left(\fint_{\mathcal{C}_{2R}}\a{\nabla\varphi}^2\right)^\frac{1}{2}.$$
Let $\delta\in(0,1)$ be arbitrary and, for each $R>0$, fix $\epsilon(R):=\delta R$.  For this choice, for each $R>0$ and $\delta\in(0,1)$,
\begin{equation}\label{sub_13}\frac{1}{R}\left(\fint^0_{-R^2}\a{(\varphi)_{t,R}-(\varphi)_R}^2\right)^{\frac{1}{2}}\lesssim \frac{1}{\delta R}\left(\fint_{\mathcal{C}_{2R}}\a{F}^2\right)^\frac{1}{2}+\delta\left(\fint_{\mathcal{C}_{2R}}\a{\nabla\varphi}^2\right)^\frac{1}{2}.\end{equation}
Since the ergodic theorem implies that, for $\en{\cdot}$-a.e. $a$,
$$\lim_{R\rightarrow\infty}\left(\fint_{\mathcal{C}_{2R}}\a{F}^2\right)^\frac{1}{2}=\en{\a{F}^2}^\frac{1}{2}\;\;\textrm{and}\;\;\lim_{R\rightarrow\infty}\left(\fint_{\mathcal{C}_{2R}}\a{\nabla\varphi}^2\right)^\frac{1}{2}=\en{\a{\nabla\varphi}^2}^\frac{1}{2},$$
it follows from (\ref{sub_13}) that, for $\en{\cdot}$-a.e. $a$, for every $\delta\in(0,1)$,
$$\limsup_{R\rightarrow\infty}\frac{1}{R}\left(\fint^0_{-R^2}\a{(\varphi)_{t,R}-(\varphi)_R}^2\right)^{\frac{1}{2}}\lesssim\delta\en{\a{\nabla\varphi}^2}^\frac{1}{2}.$$
Therefore, since $\delta\in(0,1)$ is arbitrary,
\begin{equation}\label{sub_14}\lim_{R\rightarrow 0}\frac{1}{R}\left(\fint_{\mathcal{C}_R}\a{(\varphi)_{t,R}-(\varphi)_R}^2\right)^\frac{1}{2}=\lim_{R\rightarrow 0}\frac{1}{R}\left(\fint^0_{-R^2}\a{(\varphi)_{t,R}-(\varphi)_R}^2\right)^\frac{1}{2}=0.\end{equation}
In combination, (\ref{sub_7}), (\ref{sub_8}) and (\ref{sub_14}) combine to prove (\ref{sub_5}), and thereby complete the proof of Lemma~\ref{lemma_sub_phi}.  The sublinearity of the corrector $\phi$, for $\en{\cdot}$-a.e. $a$, is then immediate from Lemma~\ref{i_lemma1}.

\subsection{The sublinearity of $\zeta$}

The fact that, for $\en{\cdot}$-a.e. $a$, for each $i,j\in\{1,\ldots,d\}$,
$$\lim_{R\rightarrow 0}\frac{1}{R^2}\left(\fint_{\mathcal{C}_R}\a{\zeta_{ij}-(\zeta_{ij})_R}^2\right)^\frac{1}{2}=\lim_{R\rightarrow\infty}\frac{1}{R^2}\left(\fint^0_{-R^2}\a{\zeta_{ij}-(\zeta_{ij})_R}^2\right)^\frac{1}{2}=0,$$
follows similarly to (\ref{new_3}).  That is, for each $i\in\{1,\ldots,d\}$, the Poincar\'e inequality and the ergodic theorem together with the Rellich-Kondrachov embedding theorem imply that the family
$$\left\{\epsilon^2\zeta_{ij}(\frac{\cdot}{\epsilon^2})-\epsilon^2\fint_{-\frac{1}{\epsilon^2}}^0\zeta_{ij}\right\}_{\epsilon\in(0,1)}$$
is compact in $L^2([0,1])$ and converges weakly to zero, as $\epsilon\rightarrow 0$, in $H^1([0,1])$.  Therefore, for $\en{\cdot}$-a.e. $a$, for each $i,j\in\{1,\ldots,d\}$, as $\epsilon\rightarrow 0$,
\begin{equation}\label{zeta_nonsense_1}\left(\epsilon^2\zeta_{ij}(\frac{\cdot}{\epsilon^2})-\epsilon^2\fint_{-\frac{1}{\epsilon^2}}^0\zeta_{ij}\right)\rightarrow 0\;\;\textrm{in}\;\;L^2([0,1]).\end{equation}

Furthermore, now exploiting the fact that $\zeta$ is a one-dimensional function, it follows from the Sobolev embedding theorem and the Arzel\`a-Ascoli theorem that, for $\en{\cdot}$-a.e. $a$, for each $i,j\in\{1,\ldots,d\}$, the family
$$\left\{\epsilon^2\zeta_{ij}(\frac{\cdot}{\epsilon^2})-\epsilon^2\fint_{-\frac{1}{\epsilon^2}}^0\zeta_{ij}\right\}_{\epsilon\in(0,1)}$$
is compact in $C^{0,\frac{1}{2}}([0,1])$ and,  by repeating the argument of (\ref{new_3}), converges weakly to zero, as $\epsilon\rightarrow 0$, in $H^1([0,1])$.  Therefore, for $\en{\cdot}$-a.e. $a$, for each $i,j\in\{1,\ldots,d\}$, as $\epsilon\rightarrow 0$,
$$\left(\epsilon^2\zeta_{ij}(\frac{\cdot}{\epsilon^2})-\epsilon^2\fint_{-\frac{1}{\epsilon^2}}^0\zeta_{ij}\right)\rightarrow 0\;\;\textrm{in}\;\;\C^{0,\frac{1}{2}}([0,1]).$$
In particular, for each $i,j\in\{1,\ldots,d\}$, as $\epsilon\rightarrow 0$,
\begin{equation}\label{zeta_nonsense_2}\a{\epsilon^2\zeta_{ij}(0)-\epsilon^2\fint_{-\frac{1}{\epsilon^2}}^0\zeta_{ij}}\rightarrow 0.\end{equation}
Hence, in combination, (\ref{zeta_nonsense_1}) and (\ref{zeta_nonsense_2}) prove after rescaling that, for $\en{\cdot}$-a.e. $a$, for each $i,j\in\{1,\ldots,d\}$,
\begin{multline*} \limsup_{R\rightarrow 0}\frac{1}{R^2}\left(\fint_{-R^2}^0\a{\zeta_{ij}-\zeta_{ij}(0)}\right)^\frac{1}{2}\leq \\ \limsup_{R\rightarrow\infty}\frac{1}{R^2}\left(\fint_{-R^2}^0\a{\zeta_{ij}-(\zeta_{ij})_R}\right)^\frac{1}{2}+\limsup_{R\rightarrow\infty}\frac{1}{R^2}\a{(\zeta_{ij})_R-\zeta_{ij}(0)}=0,\end{multline*}
which completes the argument since $\zeta$ is constant in space.

\subsection{The large-scale averages of $q$}

It is an immediately consequence of the ergodic theorem and the fact that the flux $q$ is stationary with finite energy that, for $\en{\cdot}$-a.e. $a$, for each $i\in\{1,\ldots,d\}$,
$$\lim_{R\rightarrow\infty}\left(\fint_{\mathcal{C}_R}\a{q_i}^2\right)^\frac{1}{2}=\en{\a{Q_i}^2}^\frac{1}{2},$$
which completes the argument, and the proof of Proposition~\ref{i_sublinear}.

\section{The proof of Proposition~\ref{i_excess_decay}}

The proof of Proposition~\ref{i_excess_decay} is split into five steps.  The first defines the augmented homogenization error.  The second proves that the augmented homogenization error satisfies a parabolic equation.  The third recalls some classical estimates governing the interior and boundary regularity of $\ah$-caloric functions.  The fourth uses the equation satisfied by the augmented homogenization error to derive an energy estimate.  And, finally, the fifth uses the energy estimate to complete the proof of excess decay.

In what follows, to simplify the notation, observe that it may be assumed without loss of generality that, for each $R>0$, $t\in\mathbb{R}$ and $i,j,k\in\{1,\ldots,d\}$,
$$(\phi_i)_R=(\psi_i)_{t,R}=(\sigma_{ijk})_{t,R}=\zeta_{ij}(0)=0.$$
Indeed, otherwise in the arguments to follow, at each step replace the components of the corrector, for each $R>0$, $t\in\mathbb{R}$ and $i,j,k\in\{1,\ldots,d\}$, by the normalizations defined by
\begin{equation}\label{non_normalization}\tilde{\phi}_i:=\phi_i(x,t)-(\phi_i)_R,\;\;\tilde{\psi}_i:=\psi_i-(\psi_i)_{t,R},\;\;\tilde{\sigma}_{ijk}:=\sigma_{ijk}-(\sigma_{ijk})_{t,R}\;\;\textrm{and}\;\;\tilde{\zeta}_{ij}:=\zeta_{ij}-\zeta_{ij}(0).\end{equation}
The argument now begins with the definition of the augmented homogenization error.

\subsection{The augmented homogenization error}

The analysis augmented homogenization error and its corresponding energy estimate will first be obtained on scale $R=1$.  The general results will then follow by scaling.  Suppose that $u$ is an $a$-caloric function in $\mathcal{C}_1$.  That is, in the sense of distributions, suppose that $u$ satisfies
\begin{equation}\label{w_1}u_t=\nabla\cdot a\nabla u\;\;\textrm{in}\;\;\mathcal{C}_1.\end{equation}
Then, let $\rho\in\C^\infty_c(\mathbb{R}^d)$ be a standard convolution kernel satisfying $\supp(\rho)\subset B_1$.  For each $\epsilon\in(0,\frac{1}{4})$, let $\rho^\epsilon(\cdot):=\epsilon^{-d}\rho(\frac{\cdot}{\epsilon})$ and define the spatial convolution
$$u^\epsilon(x,t)=\int_{\mathbb{R}^d}\rho^\epsilon(y-x)u(y,t)\dy\;\;\textrm{on}\;\;\mathcal{C}_{1-\epsilon}.$$
It is necessary to observe some useful energy estimates for $u$ and its convolution.

First, it is immediate from (\ref{w_1}) and the uniform ellipticity of $a$ from (\ref{i_bounded}) that
$$\norm{u_t}_{L^2([-1,0];H^{-1}(B_1))}\lesssim \int_{\mathcal{C}_1}\a{\nabla u}^2.$$
Therefore, since the convolution preserves this estimate, for each $\epsilon\in(0,\frac{1}{4})$,
\begin{equation}\label{w_time_derivative_con} \norm{u^\epsilon_t}_{L^2([-1,0];H^{-1}(B_{1-\epsilon}))}\lesssim \int_{\mathcal{C}_1}\a{\nabla u}^2.\end{equation}
It is important to keep these estimates in mind when considering the application of the constant-coefficient regularity estimates (\ref{be_1}) and (\ref{be_2}) below.

Next, although there is no convolution in time, the spatial convolution nevertheless provides some temporal regularity in the sense that, for each $\epsilon\in(0,\frac{1}{4})$,
$$u^\epsilon_t(x,t)=-\int_{\mathbb{R}^d}\nabla\rho^\epsilon(y-x)\cdot a\nabla u(y,t)\dy\;\;\textrm{in}\;\;\mathcal{C}_{1-\epsilon}.$$
Therefore, the time-derivative has a uniformly bounded energy.  That is, for each $\epsilon\in(0,\frac{1}{4})$, the Minkowski integral inequality, H\"older's inequality, the definition of the convolution kernel and the uniform ellipticity of $a$ imply that
$$\begin{aligned} \left(\int_{\mathcal{C}_{\frac{3}{4}}}\a{u^\epsilon_t}^2\right)^\frac{1}{2} \le & \left(\int_{\mathcal{C}_{1-\epsilon}}\a{\int_{\mathbb{R}^d}\nabla\rho^\epsilon(y-x)\cdot a\nabla u(y)\dy}^2\dx\dt\right)^\frac{1}{2} \\ \lesssim & \int_{\mathbb{R}^d}\left(\int_{\mathcal{C}_{1-\epsilon}}\a{\nabla\rho^\epsilon(y)}^2\a{\nabla u(y+x,t)}^2\dx\dt\right)^\frac{1}{2}\dy \\ = & \int_{\mathbb{R}^d}\a{\nabla\rho^\epsilon(y)}\left(\int_{\mathcal{C}_{1-\epsilon}}\a{\nabla u(y+x,t)}^2\dx\dt\right)^\frac{1}{2} \\ \lesssim & \frac{1}{\epsilon}\left(\int_{\mathcal{C}_1}\a{\nabla u}^2\right)^\frac{1}{2}.\end{aligned}$$

The convolution error is also well controlled by the energy of $u$.  Precisely, for each $\epsilon\in(0,\frac{1}{4})$, it follows from Jensen's inequality and the definition of the convolution kernel that
$$\begin{aligned} \left(\int_{\mathcal{C}_{\frac{3}{4}}}\a{u^\epsilon-u}^2\right)^\frac{1}{2} &=  \left(\int_{\mathcal{C}_{\frac{3}{4}}}\a{\int_0^1\int_{\mathbb{R^d}}\rho^\epsilon(y)\nabla u(x+sy,t)\cdot y \dy \ds}^2\dx\dt\right)^\frac{1}{2}\\ 
&\leq \epsilon\left(\int_{\mathbb{R}^d}\rho^\epsilon(y)\int_0^1\int_{\mathcal{C}_{\frac{3}{4}}}\a{\nabla u(x+sy,t)}^2\dx\dt\ds\dy\right)^\frac{1}{2} \\ 
&\leq   \epsilon \left(\int_{\mathcal{C}_1}\a{\nabla u}^2\right)^\frac{1}{2}.\end{aligned}$$

Lastly, it is immediate from Jensen's inequality that the convolution preserves the energy in the sense that, for each $\epsilon\in(0,\frac{1}{4})$,
$$\left(\int_{\mathcal{C}_{\frac{3}{4}}}\a{\nabla u^\epsilon}^2\right)^\frac{1}{2}\leq \left(\int_{\mathcal{C}_1}\a{\nabla u}^2\right)^\frac{1}{2}.$$

Fubini's theorem therefore implies that, for each $\epsilon\in(0,\frac{1}{4})$, there exists $r_\epsilon\in(\frac{1}{2},\frac{3}{4})$ such that
\begin{equation}\label{rep_1}\left(\int_{\partial_p\mathcal{C}_{r_\epsilon}}\a{u^\epsilon-u}^2\right)^\frac{1}{2}\lesssim \epsilon\left(\int_{\mathcal{C}_1}\a{\nabla u}^2\right)^\frac{1}{2},\end{equation}
and
\begin{equation}\label{rep_2}\left(\int_{\partial_p\mathcal{C}_{r_\epsilon}}\a{u^\epsilon_t}^2\right)^\frac{1}{2}\lesssim \frac{1}{\epsilon}\left(\int_{\mathcal{C}_1}\a{\nabla u}^2\right)^\frac{1}{2},\end{equation}
and, finally, such that
\begin{equation}\label{rep_3}\left(\int_{\partial_p\mathcal{C}_{r_\epsilon}}\a{\nabla u}^2\right)^\frac{1}{2}+\left(\int_{\partial_p\mathcal{C}_{r_\epsilon}}\a{\nabla u^\epsilon}^2\right)^\frac{1}{2}\lesssim \left(\int_{\mathcal{C}_1}\a{\nabla u}^2\right)^\frac{1}{2}.\end{equation}
It will be for this radius that the $\ah$-caloric extension of $u^\epsilon$ is constructed.

Namely, for each $\epsilon\in(0,\frac{1}{4})$, let $v^\epsilon$ denote the solution
\begin{equation}\label{w_2}\left\{\begin{array}{rll} v^\epsilon_t & =\nabla\cdot a\nabla v^\epsilon & \textrm{in}\;\;C_{r_\epsilon} \\ v^\epsilon & = u^\epsilon & \textrm{on}\;\;\partial_p C_{r_\epsilon}.\end{array}\right.\end{equation}
These functions will now come to define the augmented homogenization error after the introduction of a cutoff function.

For each $\epsilon\in(0,\frac{1}{4})$ and $\rho\in(0,\frac{1}{8})$, let $\eta^\epsilon_\rho\in\C^\infty_c(\mathbb{R}^{d+1})$ be smooth cutoff function satisfying $0\leq\eta_\rho^\epsilon\leq 1$ with, for each $x\in\mathbb{R}^d$ and $t\in\mathbb{R}$,
\begin{equation}\label{w_3}\eta_\rho^\epsilon(x,t)=\left\{\begin{array}{ll} 1 & \textrm{if}\;\;(x,t)\in\overline{\mathcal{C}}_{r_\epsilon-2\rho} \\ 0 & \textrm{if}\;\;(x,t)\in\left(\mathbb{R}^d\times(-\infty,0)\right)\setminus\mathcal{C}_{r_\epsilon-\rho}.\end{array}\right.\end{equation}
Furthermore, for each $\epsilon\in(0,\frac{1}{4})$ and $\rho\in(0,\frac{1}{8})$, for each $x\in\mathbb{R}^d$ and $t\in\mathbb{R}$,
\begin{equation}\label{w_4}\a{\nabla\eta_\rho^\epsilon(x,t)}+\a{\partial_t\eta_\rho^\epsilon(x,t)}\lesssim\frac{1}{\rho}\;\;\textrm{and}\;\;\a{\nabla^2\eta_\rho^\epsilon(x,t)}\lesssim\frac{1}{\rho^2}.\end{equation}
Then, for $\rho\in(0,\frac{1}{4})$ and $\epsilon\in(0,\frac{1}{4})$ to be specified later, define the augmented homogenization error $w$ according to the rule
\begin{equation}\label{w_5}w=u-(1+\eta_\rho^\epsilon\phi_i\partial_i)v^\epsilon.\end{equation}
The augmented homogenization error (\ref{w_5}) will now be shown to satisfy a useful parabolic equation.  The computation is motivated by the analogous computation in \cite[Lemma~2]{GNO4}, but there are significant differences owing to the parabolic setting and the use of the parabolic extended corrector $(\phi,\psi,\sigma,\zeta)$.

\subsection{The equation satisfied by the augmented homogenization error}  It is now shown that the augmented homogenization error (\ref{w_5}) satisfies
\begin{equation}\label{w_06}\left\{\begin{array}{rll} w_t-\nabla\cdot a\nabla w & = \nabla\cdot\left((1-\eta_\rho^\epsilon)(a-\ah)\nabla v^\epsilon\right)+\nabla\cdot\left((\phi_ia+\psi_i-\sigma_i)\nabla(\eta_\rho^\epsilon\partial_iv^\epsilon)\right) & \textrm{in}\;\mathcal{C}_{r_\epsilon} \\   & + \partial_t\zeta_i\cdot\nabla(\eta_\rho^\epsilon\partial_iv^\epsilon)-\phi_i (\eta_\rho^\epsilon\partial_iv^\epsilon)_t-\psi_i\Delta(\eta_\rho^\epsilon\partial_iv^\epsilon) & \\ w & = u-u^\epsilon &   \textrm{on}\;\partial_p\mathcal{C}_{r_\epsilon}.\end{array}\right.\end{equation}

Fix $\rho\in(0,\frac{1}{4})$ and $\epsilon\in(0,\frac{1}{4})$ and let $w$ be defined by (\ref{w_5}).  Since the boundary condition is immediate from the definition, it remains only to compute the equation.  First, using definition (\ref{w_5}), the gradient is defined by
$$\nabla w=\nabla u-\nabla \ve-\nabla(\phi_i\eta_\rho^\epsilon\partial_iv^\epsilon).$$
Then, because $u$ satisfies (\ref{w_1}),
\begin{equation}\label{w_6}-\nabla\cdot a\nabla w=-u_t+\nabla \cdot a\nabla v^\epsilon+\nabla\cdot a(\nabla\phi_i\eta_\rho^\epsilon\partial_iv^\epsilon+\phi_i\nabla(\eta_\rho^\epsilon\partial_iv^\epsilon)).\end{equation}

It is necessary to further analyze the term
$$\nabla \cdot a\nabla \ve+\nabla\cdot a(\nabla\phi_i\eta_\rho^\epsilon\partial_iv^\epsilon),$$
which, after adding and subtracting the unit vectors $\{e_i\}_{i\in\{1,\ldots,d\}}$, satisfies, for the fluxes $\{q_i\}_{i\in\{1,\ldots,d\}}$ defined in (\ref{i_flux}),
$$\nabla \cdot a\nabla \ve+\nabla\cdot a(\nabla\phi_i\eta_\rho^\epsilon\partial_iv^\epsilon)=\nabla\cdot ((1-\eta_\rho^\epsilon)a\nabla \ve)+\nabla\cdot (q_i \eta_\rho^\epsilon\partial_iv^\epsilon).$$
Then, after adding and subtracting the vectors $\{\ah e_i\}_{i\in\{1,\ldots,d\}}$,
\begin{multline*}\nabla \cdot a\nabla \ve+\nabla\cdot a(\nabla\phi_i\eta_\rho^\epsilon\partial_iv^\epsilon)= \\ \nabla\cdot((1-\eta_\rho^\epsilon)a\nabla \ve)+\nabla\cdot ((q_i-\ah e_i)\eta_\rho^\epsilon\partial_iv^\epsilon)+\nabla\cdot (\eta_\rho^\epsilon\ah\nabla \ve).\end{multline*}
Therefore, since $v$ satisfies (\ref{w_2}),
\begin{equation}\label{w_7}\nabla \cdot a\nabla \ve+\nabla\cdot a(\nabla\phi_i\eta_\rho^\epsilon\partial_iv^\epsilon)=\ve_t+\nabla\cdot((1-\eta_\rho^\epsilon)(a-\ah)\nabla \ve)+\nabla\cdot ((q_i-\ah e_i)\eta_\rho^\epsilon\partial_iv^\epsilon).\end{equation}
Returning to (\ref{w_6}), in view of (\ref{w_7}),
\begin{equation}\label{w_8} -\nabla\cdot a\nabla w=-u_t+\ve_t+\nabla\cdot((1-\eta_\rho^\epsilon)(a-\ah)\nabla \ve)+\nabla\cdot ((q_i-\ah e_i)\eta_\rho^\epsilon\partial_iv^\epsilon+\phi_ia\nabla(\eta_\rho^\epsilon\partial_iv^\epsilon)).\end{equation}

For the derivative in time, owing to definition (\ref{w_5}),
$$w_t=u_t-\ve_t-\phi_{i,t}\eta_\rho^\epsilon\partial_iv^\epsilon-\phi_i (\eta_\rho^\epsilon\partial_iv^\epsilon)_t,$$
which, in combination with (\ref{w_8}), yields the distributional equality
$$\begin{aligned} w_t-\nabla\cdot a\nabla w= & \nabla\cdot\left((1-\eta_\rho^\epsilon)(a-\ah)\nabla \ve+\phi_ia\nabla(\eta_\rho^\epsilon\partial_iv^\epsilon)\right)+(q_i-\ah e_i)\cdot\nabla(\eta_\rho^\epsilon\partial_iv^\epsilon) \\ + & \left((\nabla\cdot q_i)-\phi_{i,t}\right)\eta_\rho^\epsilon\partial_iv^\epsilon-\phi_i (\eta_\rho^\epsilon\partial_iv^\epsilon)_t.\end{aligned}$$
Therefore, since the correctors $\{\phi_i\}_{i\in\{1,\ldots,d\}}$ satisfy (\ref{i_corrector}), distributionally
\begin{equation}\label{w_9}\begin{aligned} w_t-\nabla\cdot a\nabla w= & \nabla\cdot\left((1-\eta_\rho^\epsilon)(a-\ah)\nabla \ve +\phi_ia\nabla(\eta_\rho^\epsilon\partial_iv^\epsilon)\right) \\ + &(q_i-\ah e_i)\cdot\nabla(\eta_\rho^\epsilon\partial_iv^\epsilon) -\phi_i (\eta_\rho^\epsilon\partial_iv^\epsilon)_t.\end{aligned}\end{equation}

It remains to analyze the term
$$(q_i-\ah e_i)\cdot\nabla(\eta_\rho^\epsilon\partial_iv^\epsilon).$$
For the correctors $\{\psi_i\}_{i\in\{1\ldots,d\}}$ satisfying (\ref{i_psi}), add and subtract the gradients $\{\nabla\psi_i\}_{i\in\{1,\ldots,d\}}$ and add and subtract the conditional expectations $\{\en{q_i\;|\;\mathcal{F}_{\mathbb{R}^d}}\}_{i\in\{1,\ldots,d\}}$ to obtain
$$\begin{aligned} (q_i-\ah e_i)\cdot\nabla(\eta_\rho^\epsilon\partial_iv^\epsilon)= & (q_i-\nabla\psi_i-\en{q_i\;|\;\mathcal{F}_{\mathbb{R}^d}})\cdot\nabla(\eta_\rho^\epsilon\partial_iv^\epsilon) \\ + & (\en{q_i\;|\;\mathcal{F}_{\mathbb{R}^d}}-\ah e_i)\cdot\nabla(\eta_\rho^\epsilon\partial_iv^\epsilon)+ \nabla\psi_i\cdot\nabla(\eta_\rho^\epsilon\partial_iv^\epsilon). \end{aligned}$$
Since the correctors $\{\sigma_i\}_{i\in\{1,\ldots,d\}}$ satisfy (\ref{i_sigma}) and the correctors $\{\zeta_i\}_{i\in\{1,\ldots,d\}}$ satisfy (\ref{i_zeta}),
\begin{equation}\label{w_10}(q_i-\ah e_i)\cdot\nabla(\eta_\rho^\epsilon\partial_iv^\epsilon)= (\nabla\cdot\sigma_i)\cdot\nabla(\eta_\rho^\epsilon\partial_iv^\epsilon)+\partial_t\zeta_i\cdot\nabla(\eta_\rho^\epsilon\partial_iv^\epsilon) + \nabla\psi_i\cdot\nabla(\eta_\rho^\epsilon\partial_iv^\epsilon).\end{equation}
Then, for each $i\in\{1,\dots,d\}$, the skew-symmetry of $\sigma_i$ proven in Lemma~\ref{i_lemma1} implies the distributional equality
\begin{equation}\label{w_11} \nabla\cdot(\sigma_i\nabla(\eta_\rho^\epsilon\partial_iv^\epsilon))=-(\nabla\cdot\sigma_i)\cdot\nabla(\eta_\rho^\epsilon\partial_iv^\epsilon).\end{equation}
Indeed, for each $i\in\{1,\ldots,d\}$, distributionally
\begin{multline*}\nabla\cdot(\sigma_i\nabla(\eta_\rho^\epsilon\partial_iv^\epsilon))=\partial_j(\sigma_{ijk}\partial_{k}(\eta_\rho^\epsilon\partial_iv^\epsilon))=\partial_j\sigma_{ijk}\partial_k(\eta_\rho^\epsilon\partial_iv^\epsilon)+\sigma_{ijk}\partial_j\partial_k(\eta_\rho^\epsilon\partial_iv^\epsilon) = \\ -\partial_j\sigma_{ikj}\partial_k(\eta_\rho^\epsilon\partial_iv^\epsilon)+\sigma_{ijk}\partial_j\partial_k(\eta_\rho^\epsilon\partial_iv^\epsilon)=-(\nabla\cdot\sigma_i)\cdot\nabla(\eta_\rho^\epsilon\partial_iv^\epsilon),\end{multline*}
where the penultimate inequality follows from the skew-symmetry of $\sigma$ and the final equality from the skew-symmetry of $\sigma$ and the equality of mixed partial derivatives.

Therefore, returning to (\ref{w_10}), the equality (\ref{w_11}) and the distributional equality
$$\nabla\psi_i\cdot\nabla(\eta_\rho^\epsilon\partial_iv^\epsilon)=\nabla\cdot(\psi_i\nabla(\eta_\rho^\epsilon\partial_iv^\epsilon))-\psi_i\Delta(\eta_\rho^\epsilon\partial_iv^\epsilon)$$
imply that
\begin{equation}\label{w_12}(q_i-\ah e_i)\cdot\nabla(\eta_\rho^\epsilon\partial_iv^\epsilon)= \nabla\cdot\left((\psi_i-\sigma_i)\nabla(\eta_\rho^\epsilon\partial_iv^\epsilon)\right)+\partial_t\zeta_i\cdot\nabla(\eta_\rho^\epsilon\partial_iv^\epsilon)-\psi_i\Delta(\eta_\rho^\epsilon\partial_iv^\epsilon).\end{equation}

Therefore, returning to (\ref{w_9}), in view of (\ref{w_12}), it follows that
\begin{equation}\label{w_13}\begin{aligned}w_t-\nabla\cdot a\nabla w= & \nabla\cdot\left((1-\eta_\rho^\epsilon)(a-\ah)\nabla \ve\right)+\nabla\cdot\left((\phi_ia+\psi_i-\sigma_i)\nabla(\eta_\rho^\epsilon\partial_iv^\epsilon)\right) \\ + &\partial_t\zeta_i\cdot\nabla(\eta_\rho^\epsilon\partial_iv^\epsilon)-\phi_i (\eta_\rho^\epsilon\partial_iv^\epsilon)_t-\psi_i\Delta(\eta_\rho^\epsilon\partial_iv^\epsilon),\end{aligned}\end{equation}
which completes the proof of (\ref{w_06}).  This equation will later be used to obtain an energy estimate for the augmented homogenization error.  However, first, it is useful to recall three classical estimates concerning the boundary and interior regularity of $\ah$-caloric functions.

\subsection{Interior and boundary estimates for $\ah$-caloric functions}  In this subsection, three classical estimates are presented to control the interior and boundary regularity of $\ah$-caloric functions.

In what follows, the boundary conditions will be assumed to be the trace of a function $\tilde{u}:\overline{\mathcal{C}}_1\rightarrow\mathbb{R}$ satisfying
\begin{equation}\label{be_assumption}\tilde{u}\in L^2([-1,0];H^1(B_1))\;\;\textrm{and}\;\;\tilde{u}_t\in L^2([-1,0];H^{-1}(B_1)).\end{equation}
The first estimate is the a priori energy estimate for the $\ah$-caloric extension of $\tilde{u}$ into $\mathcal{C}_1$.  That is, if $\tilde{v}$ satisfies
\begin{equation}\label{estimates_1}\left\{\begin{array}{rll} \tilde{v}_t & =\nabla\cdot \ah\nabla \tilde{v} & \textrm{in}\;\;\mathcal{C}_1 \\ \tilde{v} & = \tilde{u} & \textrm{on}\;\;\partial_p\mathcal{C}_1,\end{array}\right.\end{equation}
then,
\begin{equation}\label{be_1}\int_{\mathcal{C}_1}\a{\nabla\tilde{v}}^2\lesssim \int_{\mathcal{C}_1}\a{\nabla\tilde{u}}^2+\norm{\tilde{u}_t}^2_{L^2([-1,0];H^{-1}(B_1))}.\end{equation}
To prove (\ref{be_1}), let $\tilde{z}$ denote the distributional solution of
\begin{equation}\label{be_10}\left\{\begin{array}{rll} \tilde{z}_t & =\nabla\cdot \ah\nabla \tilde{z}+\nabla\cdot \ah\nabla \tilde{u}-\tilde{u}_t &  \textrm{in}\;\;\mathcal{C}_1 \\ \tilde{z} & =0 & \textrm{on}\;\;\partial\mathcal{C}_1.\end{array}\right.\end{equation}
Then, testing (\ref{be_10}) with $\tilde{z}$ and, after applying the Poincar\'e inequality, H\"older's inequality and Young's inequality and using the uniform ellipticity of $\ah$ from Lemma~\ref{i_lemma1}, it follows that
\begin{equation}\label{be_11}\int_{\mathcal{C}_1}\a{\nabla \tilde{z}}^2\lesssim \int_{\mathcal{C}_1}\a{\nabla \tilde{u}}^2+\norm{\tilde{u}_t}^2_{L^2([-1,0];H^{-1}(B_1))}.\end{equation}
However, thanks to (\ref{be_assumption}) and (\ref{estimates_1}), it is then immediate that
$$\tilde{v}=\tilde{z}+\tilde{u}.$$
Hence, with (\ref{be_11}) and the triangle inequality,
$$\int_{\mathcal{C}_1}\a{\nabla\tilde{v}}^2\lesssim \int_{\mathcal{C}_1}\a{\nabla\tilde{u}}^2+\norm{\tilde{u}_t}^2_{L^2([-1,0];H^{-1}(B_1))},$$
which proves (\ref{be_1}).

An interior regularity estimate will now be obtained for $\ah$-caloric functions.  Suppose that $\tilde{v}$ satisfies (\ref{estimates_1}) for $\tilde{u}$ satisfying (\ref{be_assumption}).  It then follows from a repeated application of the Caccioppoli inequality (\ref{i_cacc_eq}) that, for each $k\geq 0$, there exists $C(k)>0$ such that
$$\int_{B_{1-\rho}}\a{\nabla^k \tilde{v}}^2\leq \frac{C(k)}{(R\rho)^{2k}}\int_{B_1}\a{\nabla \tilde{v}}^2.$$
Therefore, by choosing $k=\frac{d}{2}$, $k=\frac{d}{2}+1$ and $k=\frac{d}{2}+2$, the Sobolev embedding theorem implies that
\begin{multline}\label{w_19}\sup_{\mathcal{C}_{1-\rho}}\left(\a{\nabla \tilde{v}}+\rho\a{\nabla^2 \tilde{v}}+\rho^2\a{\nabla^3 \tilde{v}}\right)^2\lesssim \rho^{-(d+2)}\int_{\mathcal{C}_1}\a{\nabla \tilde{v}}^2\lesssim \\ \rho^{-(d+2)}\left(\int_{\mathcal{C}_1}\a{\nabla \tilde{u}}^2+\norm{\tilde{u}_t}^2_{L^2([-1,0];H^{-1}(B_1))}\right),\end{multline}
where the final inequality follows from (\ref{be_1}).

The boundary regularity statement follows from a simplified version of Ladyzenskaja, Solonnikov and Uraltceva \cite[Theorem~9.1]{LSUbook} or, for the optimal statement, Weidemaier \cite[Theorem~3.1]{Weidemaier}.  This estimate will obtain $H^2$-regularity, and therefore requires more from the boundary condition.  In particular, this estimate explains the necessity of introducing the boundary regularization in the definition of the augmented homogenization error.  Suppose that $\tilde{u}$ satisfies the trace estimates
$$\tilde{u}\in L^2([-1,0];H^1(\partial B_1))\cap H^1(B_1\times\{-1\})\;\;\textrm{and}\;\;\tilde{u}_t\in L^2([-1,0];L^2(\partial B_1)),$$
and that $\tilde{v}$ is the $\ah$-caloric extension of $\tilde{u}$ into $\mathcal{C}_1$ in the sense of (\ref{estimates_1}).  Then, it follows from \cite[Theorem~9.1]{LSUbook} or \cite[Theorem~3.1]{Weidemaier} that
\begin{equation}\label{be_2}\int_{\mathcal{C}_1}\a{\nabla \tilde{v}}^2+\a{\nabla^2\tilde{v}}^2\lesssim \int_{\partial_p\mathcal{C}_1}\a{\nabla^\textrm{tan}\tilde{u}}^2+\int_{-1}^0\int_{\partial B_1}\a{\tilde{u}_t}^2,\end{equation}
where $\nabla^{\textrm{tan}}\tilde{u}$ denotes the tangential derivative of $\tilde{u}$ on the parabolic boundary.  In particular, $\nabla^{\textrm{tan}}\tilde{u}$ coincides with the full gradient on $B_1\times\{-1\}$.  Estimates (\ref{be_1}), (\ref{w_19}) and (\ref{be_2}) will play an important role in the energy estimate to follow.

\subsection{The energy estimate for the augmented homogenization error}  Equation (\ref{w_06}) will now be used to obtain an energy estimate for the augmented homogenization error $w$ defined in (\ref{w_5}).  Precisely, it will be shown that
\begin{equation}\begin{aligned}\label{w_130} \int_{\mathcal{C}_{r_\epsilon}}\nabla w\cdot a\nabla w\lesssim & \epsilon\int_{\mathcal{C}_1}\a{\nabla u}^2 +  \frac{\rho^\frac{2}{d}}{\epsilon^2}\int_{\mathcal{C}_1}\a{\nabla u}^2 \\ + & \frac{1}{\rho^{d+4}}\int_{\mathcal{C}_1}\left(\a{\phi}^2+\a{\psi}^2+\a{\sigma}^2\right)\int_{\mathcal{C}_1}\a{\nabla u}^2 \\ + & \left(\frac{1}{\rho^{\frac{d}{2}+3}}\left(\int_{\mathcal{C}_1}\a{\zeta}^2\right)^{\frac{1}{2}}+\frac{1}{\rho^{d+6}}\int_{\mathcal{C}_1}\a{\zeta}^2\right)\int_{\mathcal{C}_1}\a{\nabla u}^2 \\ + & \frac{1}{\rho^{d+4}}\left(\int_{\mathcal{C}_1}\a{\zeta}^2\right)^\frac{1}{2}\left(\int_{\mathcal{C}_1}\a{q}^2\right)^\frac{1}{2}\int_{\mathcal{C}_1}\a{\nabla u}^2. \end{aligned}\end{equation}
The idea is to test equation (\ref{w_13}) with $w$.  However, for this it is necessary to introduce a cutoff to ensure that $w$ vanishes along the upper boundary of the cylinder.  For each $\delta\in(0,1)$, define a smooth cutoff function $\gamma_\delta:\mathbb{R}\rightarrow\mathbb{R}$ which is non-increasing and satisfies $0\leq \gamma_\delta\leq 1$ with
\begin{equation}\label{w_014} \gamma_\delta(t)=\left\{\begin{array}{ll} 1 & \textrm{if}\;\;t\leq -\delta, \\ 0 & \textrm{if}\;\;t\geq 0.\end{array}\right.\end{equation}
Furthermore, for the Dirac mass $\delta_0$ at zero, as $\delta\rightarrow 0$,
\begin{equation}\label{w_0014} \a{(\gamma_\delta)_t}\rightharpoonup \delta_0\;\;\textrm{as distributions on $\mathbb{R}$.}\end{equation}

To begin, equation (\ref{w_13}) is tested against $\gamma_\delta w$.  Properties of the cutoff $\eta_\rho^\epsilon$ from (\ref{w_3}) and (\ref{w_4}), the uniform ellipticity of $a$ from $(\ref{i_bounded})$ and $\ah$ from Lemma~\ref{i_lemma1} and H\"older's inequality imply that, after bounding the time derivative of $\ve$ by its Hessian matrix,
\begin{equation}\begin{aligned} \label{w_14}\int_{\mathcal{C}_{r_\epsilon}}\a{(\gamma_\delta)_t}\a{w}^2+\int_{\mathcal{C}_{r_\epsilon}}\gamma_\delta\nabla w\cdot a\nabla w\lesssim & \int_{-r_\epsilon^2}^0\int_{\partial B_{r_\epsilon}}\gamma_\delta(u-u^\epsilon)\nu\cdot a(\nabla u-\nabla v^\epsilon) \\ + & \int_{B_{r_\epsilon}\times\{-{r_\epsilon}^2\}}\gamma_\delta\a{u-u^\epsilon}^2 \\ + &  \int_{\mathcal{C}_{r_\epsilon}\setminus\mathcal{C}_{r_\epsilon-2\rho}}\gamma_\delta\a{\nabla v^\epsilon}\a{\nabla w} \\ + & \sup_{\mathcal{C}_{r_\epsilon-\rho}}\left(\a{\nabla^2 v^\epsilon}+\frac{1}{R\rho}\a{\nabla v^\epsilon}\right)\int_{\mathcal{C}_{r_\epsilon}}\gamma_\delta\left(\a{\phi}+\a{\psi}+\a{\sigma}\right)\a{\nabla w} \\ + & \sup_{\mathcal{C}_{r_\epsilon-\rho}}\left(\frac{1}{\rho^2}\a{\nabla v^\epsilon}+\frac{1}{\rho}\a{\nabla^2 v^\epsilon}+\a{\nabla^3 v^\epsilon}\right) \int_{\mathcal{C}_{r_\epsilon}}\gamma_\delta\left(\a{\phi}+\a{\psi}\right)\a{w} \\ + &\a{\int_{\mathcal{C}_{r_\epsilon}}\partial_t\zeta_i\cdot\nabla(\eta_\rho^\epsilon\partial_iv^\epsilon)\gamma_\delta w},\end{aligned}\end{equation}
where $\nu$ denotes the interior normal and
$$\a{\phi}:=\left(\sum_{i=1}^d\a{\phi_i}^2\right)^\frac{1}{2},\;\; \a{\psi}:=\left(\sum_{i=1}^d\a{\psi_i}^2\right)^\frac{1}{2}\;\;\textrm{and}\;\;\a{\sigma}:=\left(\sum_{i,j,k=1}^d\a{\sigma_{ijk}}^2\right)^\frac{1}{2}.$$

For the first two boundary terms, it is immediate from the choice of $r_\epsilon\in(\frac{1}{2}, \frac{3}{4})$ in (\ref{rep_1}) and (\ref{rep_3}), the uniform ellipticity of $a$, H\"older's inequality and the estimate for the Dirichlet to Neumann map, see \cite{Brown89}, that
\begin{multline*}\int_{-r_\epsilon^2}^0\int_{\partial B_{r_\epsilon}}(u-u^\epsilon)\nu\cdot a(\nabla u-\nabla v^\epsilon)+\int_{B_{r_\epsilon}\times\{-r_\epsilon^2\}}\a{u-u^\epsilon}^2 \lesssim \\ \left(\int_{\partial_p\mathcal{C}_{r_\epsilon}}\a{u-u^\epsilon}^2\right)^\frac{1}{2}\left(\int_{\partial_p\mathcal{C}_{r_\epsilon}}\a{\nabla u}^2+\a{\nabla^\textrm{tan}u}^2\right)^\frac{1}{2}+\int_{\partial_p\mathcal{C}_{r_\epsilon}}\a{u-u^\epsilon}^2\lesssim \epsilon\int_{\mathcal{C}_1}\a{\nabla u}^2. \end{multline*}

It is then necessary to analyze the final term on the right hand side of (\ref{w_14}).  Using the definition of $w$ from (\ref{w_5}) and the fact that $\zeta$ vanishes at $t=0$ owing to (\ref{non_normalization}), it follows after integrating by parts variously in time and space that
\begin{equation}\label{w_15}\begin{aligned}\int_{\mathcal{C}_{r_\epsilon}}\partial_t\zeta_i\cdot\nabla(\eta_\rho^\epsilon\partial_iv^\epsilon)\gamma_\delta w = &\int_{\mathcal{C}_{r_\epsilon}}(\eta_\rho^\epsilon\partial_iv^\epsilon)_t\zeta_i\cdot \gamma_\delta\nabla w-\int_{\mathcal{C}_{r_\epsilon}}\zeta_i\cdot\nabla(\eta_\rho^\epsilon\partial_iv^\epsilon)(\gamma_\delta)_t w \\ - &  \int_{\mathcal{C}_{r_\epsilon}}\zeta_i\cdot\nabla(\eta_\rho^\epsilon\partial_iv^\epsilon)\gamma_\delta(u_t-v_t)-\int_{\mathcal{C}_{r_\epsilon}}\zeta_i\cdot\nabla (\eta_\rho^\epsilon\partial_iv^\epsilon)\gamma_\delta(\eta_\rho^\epsilon\phi_j\partial_j \ve)_t, \end{aligned}\end{equation}
where this equality uses the fact that the corrector $\zeta$ and the cutoff $\gamma_\delta$ are constant in space.

The first two terms of (\ref{w_15}) are bounded immediately using the definition of the cutoff $\eta_\rho^\epsilon$ from (\ref{w_3}) and (\ref{w_4}), which yields
\begin{equation}\label{w_015}\a{\int_{\mathcal{C}_{r_\epsilon}}(\eta_\rho^\epsilon\partial_iv^\epsilon)_t\zeta_i\cdot \gamma_\delta\nabla w}\lesssim \sup_{\mathcal{C}_{r_\epsilon-\rho}}\left(\frac{1}{\rho^2}\a{\nabla v^\epsilon}+\a{\nabla^3v^\epsilon}\right)\int_{\mathcal{C}_{r_\epsilon}}\gamma_\delta\a{\zeta}\a{\nabla w},\end{equation}
where
$$\a{\zeta}:=\left(\sum_{i,j=1}^d\a{\zeta_{ij}}^2\right)^\frac{1}{2}.$$
Similarly,
\begin{equation}\label{w_0015}\a{\int_{\mathcal{C}_{r_\epsilon}}(\zeta_i\cdot\nabla\left(\eta_\rho^\epsilon\partial_iv^\epsilon\right))(\gamma_\delta)_tw}\lesssim \sup_{\mathcal{C}_{r_\epsilon-\rho}}\left(\frac{1}{\rho}\a{\nabla v^\epsilon}+\a{\nabla^2 v^\epsilon}\right)\int_{\mathcal{C}_{r_\epsilon}}\a{(\gamma_\delta)_t}\a{\zeta}\a{w}.\end{equation}

It is necessary to analyze the final two terms of (\ref{w_15}).  For the first of these, using the equations (\ref{w_1}) and (\ref{w_2}) satisfied by $u$ and $v$ respectively,
$$\begin{aligned}\int_{\mathcal{C}_{r_\epsilon}}\zeta_i\cdot\nabla(\eta_\rho^\epsilon\partial_iv^\epsilon)\gamma_\delta(u_t-v_t)= & -\int_{\mathcal{C}_{r_\epsilon}}\gamma_\delta\zeta_i\cdot\left(\nabla^2(\eta_\rho^\epsilon\partial_iv^\epsilon)\cdot  a\nabla u\right) \\ +& \int_{\mathcal{C}_{r_\epsilon}}\gamma_\delta\zeta_i\cdot\left(\nabla^2(\eta_\rho^\epsilon\partial_iv^\epsilon)\cdot  \ah\nabla \ve \right).\end{aligned}$$
Therefore, using the uniform ellipticity (\ref{i_bounded}) of $a$ and the uniform ellipticity of $\ah$ from Lemma~\ref{i_lemma1}, after bounding the time derivative of $v$ by the norm of its Hessian matrix,
\begin{multline}\label{w_00015}\a{\int_{\mathcal{C}_{r_\epsilon}}\zeta_i\cdot\nabla(\eta_\rho^\epsilon\partial_iv^\epsilon)\gamma_\delta(u_t-\ve_t)}\lesssim \\ \sup_{\mathcal{C}_{r_\epsilon-\rho}}\left(\frac{1}{\rho^2}\a{\nabla v^\epsilon}+\frac{1}{\rho}\a{\nabla^2v^\epsilon}+\a{\nabla^3v^\epsilon}\right)\int_{\mathcal{C}_{r_\epsilon-\rho}}\gamma_\delta\a{\zeta}\left(\a{\nabla u}+\a{\nabla v^\epsilon}\right).\end{multline}

For the final term of (\ref{w_15}),
\begin{multline*} \int_{\mathcal{C}_{r_\epsilon}}\zeta_i\cdot\nabla (\eta_\rho^\epsilon\partial_iv^\epsilon) \gamma_\delta(\eta_\rho^\epsilon\phi_j\partial_j v^\epsilon)_t = \int_{\mathcal{C}_{r_\epsilon}}\zeta_i\cdot\nabla (\eta_\rho^\epsilon\partial_iv^\epsilon)\gamma_\delta \phi_j(\eta_\rho^\epsilon\partial_j v^\epsilon)_t \\ + \int_{\mathcal{C}_{r_\epsilon}}\zeta_i\cdot\nabla (\eta_\rho^\epsilon\partial_iv^\epsilon)\gamma_\delta\phi_{j,t}(\eta_\rho^\epsilon\partial_j v^\epsilon),\end{multline*}
and, therefore, using the equation (\ref{i_corrector}) satisfied by the correctors $\{\phi_i\}_{i\in\{1,\ldots,d\}}$,
\begin{multline*}  \int_{\mathcal{C}_{r_\epsilon}}\zeta_i\cdot\nabla (\eta_\rho^\epsilon\partial_iv^\epsilon)\gamma_\delta(\eta_\rho^\epsilon\phi_j\partial_j v^\epsilon)_t = \int_{\mathcal{C}_{r_\epsilon}}\zeta_i\cdot\nabla (\eta_\rho^\epsilon\partial_iv^\epsilon)\gamma_\delta \phi_j(\eta_\rho^\epsilon\partial_j v^\epsilon)_t \\ - \int_{\mathcal{C}_{r_\epsilon}}\gamma_\delta \nabla\left(\zeta_i\cdot \nabla (\eta_\rho^\epsilon\partial_iv^\epsilon)\eta_\rho^\epsilon\partial_jv^\epsilon\right)\cdot q_j, \end{multline*}
for the fluxes $\{q_i\}_{i\in\{1,\ldots,d\}}$ defined in (\ref{i_flux}).  Hence, after bounding the time derivative of $v$ by its Hessian matrix,
\begin{multline}\label{w_000015} \a{\int_{\mathcal{C}_{r_\epsilon}}\zeta_i\cdot\nabla (\eta_\rho^\epsilon\partial_iv^\epsilon)\gamma_\delta(\eta_\rho^\epsilon\phi_j\partial_j v^\epsilon)_t} \lesssim \\ \sup_{\mathcal{C}_{r_\epsilon-\rho}}\left(\a{\nabla v^\epsilon}\left(\frac{1}{\rho^2}\a{\nabla v^\epsilon}+\frac{1}{\rho}\a{\nabla^2v^\epsilon}+\a{\nabla^3v^\epsilon}\right)+\left(\frac{1}{\rho}\a{\nabla v^\epsilon}+\a{\nabla^2v^\epsilon}\right)^2\right)\int_{\mathcal{C}_{r_\epsilon}}\gamma_\delta\a{\zeta}\a{q} \\ +\sup_{\mathcal{C}_{r_\epsilon-\rho}}\left(\left(\frac{1}{\rho}\a{\nabla v^\epsilon}+\a{\nabla^2v^\epsilon}\right)\left(\frac{1}{\rho^2}\a{\nabla v^\epsilon}+\a{\nabla^3v^\epsilon}\right)\right)\int_{\mathcal{C}_{r_\epsilon}}\gamma_\delta\a{\zeta}\a{\phi},\end{multline}
where
$$\a{q}:=\left(\sum_{i=1}^d\a{q_i}^2\right)^\frac{1}{2}.$$

Therefore, in view of (\ref{w_time_derivative_con}) and (\ref{w_14}), it follows from the uniform ellipticity of $a$, the definition of $\gamma_\delta$, the Poincar\'e inequality in space, H\"older's inequality and Young's inequality that
\begin{equation}\begin{aligned} \label{w_21}  \int_{\mathcal{C}_{r_\epsilon}}\a{(\gamma_\delta)_t}w^2+\int_{\mathcal{C}_{r_\epsilon}}\gamma_\delta\nabla w\cdot a\nabla w\lesssim & \epsilon\int_{\mathcal{C}_1}\a{\nabla u}^2 \\ + & \int_{\mathcal{C}_{r_\epsilon}\setminus \mathcal{C}_{r_\epsilon-2\rho}}\a{\nabla v^\epsilon}^2 \\ + & \frac{1}{\rho^{d+4}}\int_{\mathcal{C}_1}\left(\a{\phi}^2+\a{\psi}^2+\a{\sigma}^2\right)\int_{\mathcal{C}_1}\a{\nabla u}^2 \\ + &\a{\int_{\mathcal{C}_{r_\epsilon}}\partial_t\zeta_i\cdot\nabla(\eta_\rho^\epsilon\partial_iv^\epsilon)\gamma_\delta w}. \end{aligned}\end{equation}
For the second term of (\ref{w_21}), the choice of the radius $r_\epsilon\in(\frac{1}{2}, \frac{3}{4})$ satisfying (\ref{rep_2}) and estimate (\ref{be_1}) for $\ah$-caloric functions imply that
$$\int_{\mathcal{C}_{r_\epsilon}}\a{\nabla v^\epsilon}^2+\a{\nabla^2v^\epsilon}^2\lesssim \int_{\partial_p\mathcal{C}_{r_\epsilon}}\a{\nabla^{\textrm{tan}}u^\epsilon}^2+\int_{-r_\epsilon^2}^0\int_{\partial B_{r_\epsilon}}\a{u^\epsilon_t}^2\lesssim \frac{1}{\epsilon^2}\int_{\mathcal{C}_{r_\epsilon}}\a{\nabla u^\epsilon}^2\leq\frac{1}{\epsilon^2}\int_{\mathcal{C}_1}\a{\nabla u}^2.$$
Therefore, it follows from the Sobolev embedding theorem that
$$\int_{\mathcal{C}_{r_\epsilon}\setminus \mathcal{C}_{r_\epsilon-2\rho}}\a{\nabla v^\epsilon}^2\lesssim \left(\int_{\mathcal{C}_{r_\epsilon}}\chi_{\mathcal{C}_{r_\epsilon}\setminus \mathcal{C}_{r_\epsilon-2\rho}}\right)^\frac{2}{d}\left(\int_{\mathcal{C}_{r_\epsilon}}\a{\nabla v^\epsilon}^\frac{2d}{d-2}\right)^\frac{d-2}{d}\leq \frac{\rho^\frac{2}{d}}{\epsilon^2}\int_{\mathcal{C}_1}\a{\nabla u}^2,$$
where $\chi_{\mathcal{C}_{r_\epsilon}\setminus \mathcal{C}_{r_\epsilon-2\rho}}$ is the indicator function of the set $(\mathcal{C}_{r_\epsilon}\setminus \mathcal{C}_{r_\epsilon-2\rho})$, and where the argument is only written for the case $d\geq 3$, since the modifications necessary for the cases $d=1$ and $d=2$ are straightforward and rely only upon the Sobolev embedding theorem.

For the final term of (\ref{w_21}), estimates (\ref{w_015}), (\ref{w_0015}), (\ref{w_00015}) and (\ref{w_000015}) together with estimates (\ref{w_time_derivative_con}) and  (\ref{w_19}), where H\"older's inequality is used for the final term, prove that, since removing $\gamma_\delta$ from the final three terms of the right hand side increases their magnitude,
\begin{equation}\begin{aligned}\label{w_22} \a{\int_{\mathcal{C}_{r_\epsilon}}\partial_t\zeta_i\cdot\nabla(\eta_\rho^\epsilon\partial_iv^\epsilon)\gamma_\delta w}\lesssim & \frac{1}{\rho^{\frac{d}{2}+3}}\int_{\mathcal{C}_{r_\epsilon}}\gamma_\delta\a{\zeta}\a{\nabla w}\left(\int_{\mathcal{C}_1}\a{\nabla u}^2\right)^\frac{1}{2} \\ + & \frac{1}{\rho^{\frac{d}{2}+1}}\int_{\mathcal{C}_{r_\epsilon}}\a{(\gamma_\delta)_t}\a{\zeta}\a{w}\left(\int_{\mathcal{C}_1}\a{\nabla u}^2\right)^\frac{1}{2} \\ + & \frac{1}{\rho^{\frac{d}{2}+3}}\left(\int_{\mathcal{C}_1}\a{\zeta}^2\right)^\frac{1}{2}\int_{\mathcal{C}_1}\a{\nabla u}^2 \\ + & \frac{1}{\rho^{d+4}}\left(\int_{\mathcal{C}_1}\a{\zeta}^2\right)^\frac{1}{2}\left(\int_{\mathcal{C}_1}\a{q}^2\right)^\frac{1}{2}\int_{\mathcal{C}_1}\a{\nabla u}^2 \\ + & \frac{1}{\rho^{d+5}}\left(\int_{\mathcal{C}_1}\a{\zeta}^2\right)^\frac{1}{2}\left(\int_{\mathcal{C}_1}\a{\phi}^2\right)^\frac{1}{2}\int_{\mathcal{C}_1}\a{\nabla u}^2.\end{aligned}\end{equation}
Therefore, following an application of H\"older's inequality and then Young's inequality, it follows from (\ref{w_21}) and (\ref{w_22}) that
\begin{equation}\begin{aligned}\label{w_23} \int_{\mathcal{C}_{r_\epsilon}}\gamma_\delta\nabla w\cdot a\nabla w\lesssim &  \epsilon\int_{\mathcal{C}_1}\a{\nabla u}^2 +  \frac{\rho^\frac{2}{d}}{\epsilon^2}\int_{\mathcal{C}_1}\a{\nabla u}^2 \\  + & \frac{1}{\rho^{d+4}}\int_{\mathcal{C}_1}\left(\a{\phi}^2+\a{\psi}^2+\a{\sigma}^2\right)\int_{\mathcal{C}_1}\a{\nabla u}^2 \\ + & \left(\frac{1}{\rho^{\frac{d}{2}+3}}\left(\int_{\mathcal{C}_1}\a{\zeta}^2\right)^{\frac{1}{2}}+\frac{1}{\rho^{d+6}}\int_{\mathcal{C}_1}\a{\zeta}^2\right)\int_{\mathcal{C}_1}\a{\nabla u}^2 \\ + & \frac{1}{\rho^{d+4}}\left(\int_{\mathcal{C}_1}\a{\zeta}^2\right)^\frac{1}{2}\left(\int_{\mathcal{C}_1}\a{q}^2\right)^\frac{1}{2}\int_{\mathcal{C}_1}\a{\nabla u}^2 \\ + & \frac{1}{\rho^{d+2}}\int_{\mathcal{C}_1}\a{(\gamma_\delta)_t}\a{\zeta}^2\int_{\mathcal{C}_1}\a{\nabla u}^2. \end{aligned}\end{equation}

In view of the construction of $\gamma_\delta$ from (\ref{w_014}), and owing to the distributional convergence (\ref{w_0014}), the fact that $\zeta$ vanishes at $t=0$ thanks to (\ref{non_normalization}) implies, for $\en{\cdot}$-a.e. $a$,
$$\lim_{\delta\rightarrow 0}\int_{\mathcal{C}_1}\a{(\gamma_\delta)_t}\a{\zeta}^2=0.$$
Therefore, after passing to the limit $\delta\rightarrow 0$ in (\ref{w_23}), the construction of $\gamma_\delta$ in (\ref{w_014}) implies that
\begin{equation}\begin{aligned}\label{w_24} \int_{\mathcal{C}_{r_\epsilon}}\nabla w\cdot a\nabla w\lesssim & \epsilon\int_{\mathcal{C}_1}\a{\nabla u}^2 +  \frac{\rho^\frac{2}{d}}{\epsilon^2}\int_{\mathcal{C}_1}\a{\nabla u}^2 \\ + & \frac{1}{\rho^{d+4}}\int_{\mathcal{C}_1}\left(\a{\phi}^2+\a{\psi}^2+\a{\sigma}^2\right)\int_{\mathcal{C}_1}\a{\nabla u}^2 \\ + & \left(\frac{1}{\rho^{\frac{d}{2}+3}}\left(\int_{\mathcal{C}_1}\a{\zeta}^2\right)^{\frac{1}{2}}+\frac{1}{\rho^{d+6}}\int_{\mathcal{C}_1}\a{\zeta}^2\right)\int_{\mathcal{C}_1}\a{\nabla u}^2 \\ + & \frac{1}{\rho^{d+4}}\left(\int_{\mathcal{C}_1}\a{\zeta}^2\right)^\frac{1}{2}\left(\int_{\mathcal{C}_1}\a{q}^2\right)^\frac{1}{2}\int_{\mathcal{C}_1}\a{\nabla u}^2, \end{aligned}\end{equation}
which completes the proof of (\ref{w_130}).

To obtain (\ref{w_24}) for an arbitrary radius $R>0$, suppose that $u$ is an $a$-caloric function on $B_R$.  Then, for each $\epsilon\in(0,\frac{R}{4})$, there exists a radius $R_\epsilon\in (\frac{R}{2},\frac{3R}{4})$ and a cutoff function $\eta^{R_\epsilon}_\rho$ with $0\leq \eta^{R_\epsilon}_\rho\leq 1$ and such that
\begin{equation}\label{w_cutoff_rescaled}\eta^{R_\epsilon}_\rho(x,t)=\left\{\begin{array}{ll} 1 & \textrm{if}\;\;(x,t)\in \overline{\mathcal{C}}_{R_\epsilon-2\rho R_\epsilon} \\ 0 & \textrm{if}\;\;(x,t)\in(\mathbb{R}^d\times(-\infty,0))\setminus\mathcal{C}_{R_\epsilon-\rho R_\epsilon},\end{array}\right.\end{equation}
which, for the $\ah$-caloric extension $v^\epsilon$ of $u^\epsilon$ into $\mathcal{C}_{R_\epsilon}$, define the corresponding augmented homogenization error
$$w=u-(1+\eta^{R_\epsilon}_\rho\phi_i\partial_i)v^\epsilon.$$
Then, for each $\epsilon\in(0,\frac{R}{4})$, after performing the rescalings
$$(\tilde{w},\tilde{u},\tilde{v}^\epsilon,\tilde{q})(\cdot,\cdot)=\frac{1}{R_\epsilon}(w,u,v^\epsilon,q)(R_\epsilon\cdot,R_\epsilon^2\cdot),$$
it follows using equations (\ref{i_corrector}), (\ref{i_psi}), (\ref{i_sigma}) and (\ref{i_zeta}) that the correctors rescale according to the rules
$$(\tilde{\phi}, \tilde{\psi}, \tilde{\sigma})=\frac{1}{R_\epsilon}(\phi, \psi, \sigma)(R_\epsilon\cdot, R_\epsilon^2\cdot)\;\;\textrm{and}\;\;\tilde{\zeta}(\cdot)=\frac{1}{R_\epsilon^2}\zeta(R^2_\epsilon\cdot).$$

Hence, after applying (\ref{w_24}) and returning to the original scaling, it follows that, for each $\epsilon\in(0,\frac{R}{4})$ and $\rho\in(0,\frac{1}{8})$,
\begin{equation}\begin{aligned}\label{w_end} \fint_{\mathcal{C}_{R_\epsilon}}\nabla w\cdot a\nabla w\lesssim & \epsilon\fint_{\mathcal{C}_R}\a{\nabla u}^2 +  \frac{\rho^\frac{2}{d}}{\epsilon^2}\fint_{\mathcal{C}_R}\a{\nabla u}^2 \\ + & \frac{1}{R^2\rho^{d+4}}\fint_{\mathcal{C}_R}\left(\a{\phi}^2+\a{\psi}^2+\a{\sigma}^2\right)\fint_{\mathcal{C}_R}\a{\nabla u}^2 \\ + & \left(\frac{1}{R^2\rho^{\frac{d}{2}+3}}\left(\int_{\mathcal{C}_R}\a{\zeta}^2\right)^{\frac{1}{2}}+\frac{1}{R^4\rho^{d+6}}\int_{\mathcal{C}_R}\a{\zeta}^2\right)\int_{\mathcal{C}_R}\a{\nabla u}^2 \\ + & \frac{1}{R^2\rho^{d+4}}\left(\fint_{\mathcal{C}_R}\a{\zeta}^2\right)^\frac{1}{2}\left(\fint_{\mathcal{C}_R}\a{q}^2\right)^\frac{1}{2}\fint_{\mathcal{C}_R}\a{\nabla u}^2, \end{aligned}\end{equation}
which is the general form of the energy estimate that will be used in the proof of excess decay to follow.

\subsection{The proof of excess decay}

The energy estimate will now be used to prove the excess decay of Proposition~\ref{i_excess_decay}.  Fix $R>0$ and suppose that $u$ is an $a$-caloric function $\mathcal{C}_R$.  Then, for each $\epsilon\in(0,\frac{R}{4})$ and $\rho\in(0,\frac{1}{8})$, choose a radius $R_\epsilon\in(\frac{R}{2},\frac{3R}{4})$ and a cutoff $\eta^{R_\epsilon}_\rho$ such that, for the $\ah$-caloric extension $v^\epsilon$ of $u^\epsilon$ into $\mathcal{C}_{R_\epsilon}$, the conclusion of (\ref{w_end}) is satisfied for the augmented homogenization error $w$ defined by
\begin{equation}\label{e_1}w=u-(1+\eta^{R_\epsilon}_\rho\phi_i\partial_i)v^\epsilon\;\;\textrm{in}\;\;\mathcal{C}_{R_\epsilon}.\end{equation}
The proof of excess decay will now proceed in four steps.

{\textbf{Step 1.}}  In the first step of the proof, it will be shown that, for any $\delta>0$, there exists $C_2=C_2(d,\lambda,\delta)>0$ such that, whenever, for each $i\in\{1,\ldots,d\}$,
\begin{equation}\label{e_2}\left(\fint_{\mathcal{C}_{R}}\a{q_i}^2\right)^\frac{1}{2}\leq 2\en{\a{q_i}^2}^\frac{1}{2},\end{equation}
and
\begin{equation}\label{e_3}\frac{1}{R}\left(\fint_{\mathcal{C}_{R}}\a{\phi_i}^2+\a{\psi_i}^2+\a{\sigma_i}^2\right)^\frac{1}{2}+\frac{1}{R^2}\left(\fint_{\mathcal{C}_{R}}\a{\zeta_i}^2\right)^\frac{1}{2}\leq\frac{1}{C_2},\end{equation}
then
\begin{equation}\label{e_4}\fint_{\mathcal{C}_{R_\epsilon}}\nabla w\cdot a\nabla w\lesssim \delta\fint_{\mathcal{C}_R}\nabla u\cdot a\nabla u.\end{equation}
The proof is a simple consequence of estimate (\ref{w_end}) and the definition (\ref{e_1}).

Fix $\delta>0$.  Then, assuming that (\ref{e_2}) and (\ref{e_3}) are satisfied for some $C_2>0$ to be fixed later, estimate (\ref{w_24}), the choice (\ref{e_1}) and the uniform ellipticity of $a$ imply that
\begin{multline*}\fint_{\mathcal{C}_{R_1}}\nabla w\cdot a\nabla w\lesssim \\ \left(\epsilon+\frac{\rho^\frac{2}{d}}{\epsilon^2}+\frac{1}{C_2^2\rho^{d+4}}+\frac{1}{C_2\rho^{\frac{d}{2}+3}}+\frac{1}{C_2^2\rho^{d+6}}+\frac{1}{C_2\rho^{d+4}}\en{\a{q}^2}^\frac{1}{2} \right)\fint_{\mathcal{C}_{R_1}}\nabla u\cdot a\nabla u.\end{multline*}
Therefore, first choose $\epsilon_0\in(0,\frac{1}{4})$ satisfying
$$\epsilon_0<\frac{1}{3}\delta.$$
Then, choose $\rho_0\in(0,\frac{1}{8})$ sufficiently small so as to guarantee that
$$\frac{\rho^\frac{2}{d}}{\epsilon_0^2}<\frac{1}{3}\delta.$$
Finally, fix $C_2>0$ large enough to ensure that
$$\left(\frac{1}{C_2^2\rho_0^{d+4}}+\frac{1}{C_2\rho_0^{\frac{d}{2}+3}}+\frac{1}{C_2^2\rho_0^{d+6}}+\frac{1}{C_2\rho_0^{d+4}}\en{\a{q^2}}^\frac{1}{2} \right)<\frac{\delta}{3}.$$
Then, it follows that, for this choice of $\epsilon_0$, $\rho_0$ and $C_2$,
$$\fint_{\mathcal{C}_{R_{\epsilon_0}}}\nabla w\cdot a\nabla w\lesssim \delta\fint_{\mathcal{C}_{R}}\nabla u\cdot a\nabla u,$$
which proves (\ref{e_4}).

In particular, since $\rho_0\in(0,\frac{1}{8})$ and $R_{\epsilon_0}\in(\frac{R}{2},\frac{3R}{4})$, using the definition (\ref{w_cutoff_rescaled}) of the cutoff $\eta^{R_{\epsilon_0}}_{\rho_0}$, it follows that
\begin{equation}\label{e_5}w=u-(1+\phi_i\partial_i)v^{\epsilon_0}\;\;\textrm{on}\;\; \mathcal{C}_\frac{R}{4}.\end{equation}
Therefore, since $R_{\epsilon_0}\in(\frac{R}{2},\frac{3R}{4})$, it follows from (\ref{e_4}) and (\ref{e_5}) that, for any $\delta>0$ there exists $C_2=C_2(d,\lambda)>0$ such that, whenever (\ref{e_3}) and (\ref{e_4}) are satisfied, for each $r\in(0,\frac{R}{4}]$,
\begin{equation}\label{e_6}\fint_{\mathcal{C}_r}(\nabla u-\nabla v^{\epsilon_0}-\nabla(\phi_i\partial_i v^{\epsilon_0}))\cdot a(\nabla u-\nabla v^{\epsilon_0}-\nabla(\phi_i\partial_i v^{\epsilon_0}))\lesssim\delta\left(\frac{R}{r}\right)^{d+2}\fint_{\mathcal{C}_{R}}\nabla u\cdot a\nabla u.\end{equation}
This completes the first step of the proof.

{\textbf{Step 2.}}  The second step will show that the left hand side of (\ref{e_6}) is a good approximation for the excess by using the interior regularity of $\ah$-caloric functions and the Caccioppoli inequality.  To simplify the notation in what follows, define
$$v:=v^{\epsilon_0}\;\;\textrm{in}\;\;\mathcal{C}_{R_{\epsilon_0}}.$$
Then, form the decomposition
\begin{equation}\label{e_7}\nabla u-\nabla v-\nabla(\phi_i\partial_iv)=\nabla u-\nabla v(0,0)(I_d+\nabla\phi)+\left(\nabla v(0,0)-\nabla v\right)(I_d+\nabla\phi)-\phi_i\nabla(\partial_iv),\end{equation}
where $I_d$ denotes the $(d\times d)$-identity matrix and, for each $i,j\in\{1,\ldots,d\}$,
$$\left(\nabla\phi\right)_{ij}:=\partial_j\phi_i.$$
After fixing $\xi_0=\nabla v(0,0)$, use (\ref{e_7}), the triangle inequality and Young's inequality to prove that, in $\mathcal{C}_r$ for any $r\in(0,\frac{R}{4}]$,
\begin{equation}\label{e_8}\a{\nabla u-\xi_0-\nabla\phi_{\xi_0}}^2\lesssim \a{\nabla w}^2+\sup_{\mathcal{C}_r}\left(\a{\nabla v-\nabla v(0,0)}\right)^2\a{I+\nabla\phi}^2+\sup_{\mathcal{C}_r}\left(\a{\nabla(\partial_iv)}\right)^2\a{\phi_i}^2.\end{equation}

Estimate (\ref{w_19}) implies that, after bounding the time derivative of $v$ by the norm of its Hessian matrix, and using the uniform ellipticity of $a$ and the choice $R_{\epsilon_0}\in(\frac{R}{2}, \frac{3R}{4})$, for each $r\in(0,\frac{R}{4}]$,
\begin{multline}\label{e_9}\sup_{\mathcal{C}_r}\left(\a{\nabla v-\nabla v(0,0)}\right)\lesssim \sup_{\mathcal{C}_r}\left(r\a{\nabla^2v}+r^2\a{\nabla^3v}\right)^2 \\ \lesssim \left(\frac{r}{R}\right)^2\sup_{\mathcal{C}_\frac{R}{4}}\left(\frac{R}{4}\a{\nabla^2v}+\left(\frac{R}{4}\right)^2\a{\nabla^3v}\right)^2\lesssim \left(\frac{r}{R}\right)^2\fint_{\mathcal{C}_R}\nabla u \cdot a\nabla u.\end{multline}
Similarly, for each $i\in\{1,\ldots,d\}$, using estimate (\ref{w_19}), the uniform ellipticity of $a$ and $R_{\epsilon_0}\in(\frac{R}{2},\frac{3R}{4})$, it follows that, for each $r\in(0,\frac{R}{4}]$,
\begin{equation}\label{e_10}\sup_{\mathcal{C}_r}\left(\a{\nabla(\partial_iv)}\right)^2\lesssim\sup_{\mathcal{C}_\frac{R}{4}}\left(\a{\nabla(\partial_iv)}\right)^2 \lesssim \frac{1}{R^2}\fint_{\mathcal{C}_R}\nabla u \cdot a\nabla u.\end{equation}

Finally, since for each $i\in\{1,\ldots,d\}$ the $a$-caloric coordinate  $(x_i+\phi_i)$ satisfies
$$\partial_t(x_i+\phi_i)=\nabla\cdot a\nabla(x_i+\phi_i)\;\;\textrm{in}\;\;\mathbb{R}^{d+1},$$
the Caccioppoli inequality (\ref{i_cacc_eq}) implies that, for each $r\in(0,\frac{R}{4}]$,
\begin{equation}\label{e_11}\fint_{\mathcal{C}_r}\a{e_i+\nabla\phi_i}^2\lesssim\frac{1}{(2r)^2}\fint_{\mathcal{C}_{2r}}\a{\phi_i}^2+\frac{1}{(2r)^2}\fint_{\mathcal{C}_{2r}}x_i^2\lesssim\frac{1}{(2r)^2}\fint_{\mathcal{C}_{2r}}\a{\phi_i}^2+1.\end{equation}
Therefore, returning to (\ref{e_8}), estimates (\ref{e_9}), (\ref{e_10}) and (\ref{e_11}), with the uniform ellipticity of $a$ and the choice $R_{\epsilon_0}\in(\frac{R}{2}, \frac{3R}{4})$, imply that, for each $r\in(0,\frac{R}{4}]$,
\begin{equation}\begin{aligned}\label{e_12}\fint_{\mathcal{C}_r}(\nabla u-\xi_0-\nabla\phi_{\xi_0})\cdot a (\nabla u-\xi_0-\nabla\phi_{\xi_0})\lesssim & \left(\frac{R}{r}\right)^{d+2} \fint_{\mathcal{C}_{R_{\epsilon_0}}}\nabla w\cdot a\nabla w \\ +& \left(\frac{r}{R}\right)^2\left(\frac{1}{(2r)^2}\fint_{\mathcal{C}_{2r}}\a{\phi_i}^2+1\right)\fint_{\mathcal{C}_R}\nabla u \cdot a\nabla u.\end{aligned}\end{equation}
This completes the second step.

{\textbf{Step 3.}}  In the third step, inequality (\ref{e_12}) will be combined with (\ref{e_4}) to prove the excess decay along a subsequence.  Namely, for every $\alpha\in(0,1)$, it will be shown that there exists $C_0=C_0(d,\lambda,\alpha)>0$ and $\theta_0=\theta_0(\alpha,d,\lambda)\in(0,\frac{1}{4})$ such that, if $r_1=\theta_0R$ and if, for each $r\in[r_1,R]$,
\begin{equation}\label{e_13}\left(\fint_{\mathcal{C}_r}\a{q_i}^2\right)^\frac{1}{2}\leq 2\en{\a{q_i}^2}^\frac{1}{2},\end{equation}
and
\begin{equation}\label{e_14}\frac{1}{r}\left(\fint_{\mathcal{C}_r}\a{\phi_i}^2+\a{\psi_i}^2+\a{\sigma_i}^2\right)^\frac{1}{2}+\frac{1}{r^2}\left(\fint_{\mathcal{C}_r}\a{\zeta_i}^2\right)^\frac{1}{2}\leq\frac{1}{C_0},\end{equation}
then
\begin{equation}\label{e_15}\textrm{Exc}\left(u;r_1\right)\leq\left(\frac{r_1}{R}\right)^{2\alpha}\textrm{Exc}(u;R).\end{equation}
Notice that the inequality appearing in (\ref{e_15}) is exact.

Let $\delta>0$ be arbitrary.  In view of (\ref{e_4}), there exists $C_2=C_2(\delta,d,\lambda)\geq 1$ such that, whenever (\ref{e_13}) and (\ref{e_14}) are satisfied for the constant $C_2$, then, since $R_{\epsilon_0}\in(\frac{R}{2}, \frac{3R}{4})$,
$$\fint_{\mathcal{C}_{R_{\epsilon_0}}}\nabla w\cdot a\nabla w\lesssim \delta\fint_{\mathcal{C}_R}\nabla u\cdot a\nabla u.$$
Therefore, it follows from inequality (\ref{e_12}) and (\ref{e_14}) that, for $C_3=C_3(d,\lambda)>0$, since $R_{\epsilon_0}\in(\frac{R}{2}, \frac{3R}{4})$,
\begin{equation}\label{e_16}\fint_{\mathcal{C}_{r_1}}(\nabla u-\xi_0-\nabla\phi_{\xi_0})\cdot a (\nabla u-\xi_0-\nabla\phi_{\xi_0})\leq C_3\left(\delta\left(\frac{R}{r_1}\right)^{d+2}+\left(\frac{r_1}{R}\right)^2\right)\fint_{\mathcal{C}_R}\nabla u \cdot a\nabla u.\end{equation}
Choose $\theta_0\in(0,\frac{1}{4})$ sufficiently small so as to guarantee
\begin{equation}\label{e_17}C_3\theta_0^2\leq \frac{1}{2}\theta_0^{2\alpha},\end{equation}
which is possible because $\alpha\in(0,1)$, and choose $\delta_0>0$ sufficiently small so as to guarantee
\begin{equation}\label{e_18}C_3\delta_0\theta_0^{-(d+2)}\leq \frac{1}{2}\theta_0^{2\alpha}.\end{equation}

It is then immediate from (\ref{e_16}) that, by choosing $\theta_0$ as in (\ref{e_17}) and choosing $C_0:=C_2(\delta_0,d,\lambda)$ for $\delta_0$ defined in (\ref{e_18}), whenever (\ref{e_13}) and (\ref{e_14}) are satisfied for the constant $C_0$ and $r_1=\theta_0R$,
\begin{equation}\label{e_19}\fint_{\mathcal{C}_{r_1}}(\nabla u-\xi_0-\nabla\phi_{\xi_0})\cdot a (\nabla u-\xi_0-\nabla\phi_{\xi_0})\leq \theta_0^{2\alpha}\fint_{\mathcal{C}_R}\nabla u \cdot a\nabla u=\left(\frac{r_1}{R}\right)^{2\alpha}\fint_{\mathcal{C}_R}\nabla u \cdot a\nabla u.\end{equation}
Since the excess is defined, for each $R>0$, by
$$\textrm{Exc}(u;R)=\inf_{\xi\in\mathbb{R}^d}\fint_{\mathcal{C}_R}(\nabla u-\xi-\nabla\phi_\xi)\cdot a(\nabla u-\xi-\nabla\phi_\xi),$$
inequality (\ref{e_19}) implies that
\begin{equation}\label{e_20}\textrm{Exc}(u;r_1)\leq \left(\frac{r_1}{R}\right)^{2\alpha}\fint_{\mathcal{C}_R}\nabla u \cdot a\nabla u.\end{equation}
However, because the left hand side of inequality (\ref{e_20}) is invariant with respect to the addition of an arbitrary $a$-caloric gradient $(\xi+\nabla\phi_\xi)$ in the sense that, with (\ref{e_20}), for every $\xi\in\mathbb{R}^d$,
\begin{equation}\label{e_21}\textrm{Exc}(u;r_1)=\textrm{Exc}(u-(\xi\cdot x+\phi_\xi);r_1)\leq\left(\frac{r_1}{R}\right)^{2\alpha}\fint_{\mathcal{C}_R}(\nabla u-\xi-\nabla\phi_\xi) \cdot a\nabla (u-\xi-\nabla\phi_\xi),\end{equation}
taking an infimum on the right hand side with respect to $\xi\in\mathbb{R}^d$ yields
\begin{equation}\label{e_21}\textrm{Exc}(u;r_1)\leq\left(\frac{r_1}{R}\right)^{2\alpha}\textrm{Exc}(u;R),\end{equation}
which completes the proof of (\ref{e_15}), and the argument's third step.

{\textbf{Step 4.}}  The final step completes the proof using (\ref{e_21}) and an iteration argument.  Fix $r_1<R$ such that, for $C_0>0$ defined following (\ref{e_18}), both (\ref{e_13}) and (\ref{e_14}) are satisfied for the constant $C_0$ for every $r\in[r_1,R]$.  It will be shown that, in this case,
\begin{equation}\label{e_22}\textrm{Exc}(u;r)\lesssim \left(\frac{r}{R}\right)^{2\alpha}\textrm{Exc}(u;R).\end{equation}
Fix $\theta_0$ as defined in (\ref{e_14}).  If $r\geq R\theta_0$, then using the definition of the excess, for $C=C(\theta_0)>0$,
\begin{multline}\label{e_23}\mathrm{Exc}(u;r)\leq \left(\frac{R}{r}\right)^d\mathrm{Exc}(u;R)=\left(\frac{R}{r}\right)^{d+2\alpha}\left(\frac{r}{R}\right)^{2\alpha}\mathrm{Exc}(u;R) \\ \leq \theta_0^{-(d+2\alpha)}\left(\frac{r}{R}\right)^{2\alpha}\mathrm{Exc}(u;R) \leq C\left(\frac{r}{R}\right)^{2\alpha}\mathrm{Exc}(u;R).\end{multline}
If $r<\theta_0R$, then let $n$ be the unique positive integer satisfying $\theta_0^{n-1}R \leq r < \theta_0^nR$.  Proceeding inductively, and relying upon the fact that (\ref{e_21}) obtains an exact inequality, for constants $C=C(\theta_0)>0$ which can change between inequalities,
\begin{multline}\label{e_24}
 \mathrm{Exc}(u;r) \leq C \mathrm{Exc}(u;\theta_0^n R) \leq C (\theta_0^{n})^{2\alpha} \mathrm{Exc}(u;R) = \\ C\theta_0^{2\alpha} (\theta_0^{n-1})^{2\alpha} \mathrm{Exc}(u;R) \leq  C\left( \frac{r}{R} \right)^{2\alpha} \mathrm{Exc}(u;R). \end{multline}
In combination, (\ref{e_23}) and (\ref{e_24}) prove (\ref{e_22}) and complete the proof of Proposition~\ref{i_excess_decay}.

\section{The proof of Lemma~\ref{i_cacc}}

Fix a coefficient field $a$ satisfying (\ref{i_bounded}).  Fix $R>0$ and suppose that $u$ is a distributional solution of
\begin{equation}\label{cacc_eq} u_t=\nabla\cdot a\nabla u\;\;\textrm{in}\;\;\mathcal{C}_R.\end{equation}
Let $c\in\mathbb{R}$ and $\rho\in(0,\frac{R}{2})$ be arbitrary.  The Caccioppoli inequality is obtained by testing equation (\ref{cacc_eq}) with $\eta^2(u-c)$ for an appropriately chosen cutoff function $\eta$.

Precisely, fix $\eta\in\C^\infty_c(\mathbb{R}^{d+1})$ satisfying $0\leq \eta\leq 1$ and, for $x\in\mathbb{R}^d$ and $t\in\mathbb{R}$,
$$\eta(x,t)=\left\{\begin{array}{ll} 1 & \textrm{if}\;\;(x,t)\in \overline{B}_{R-\rho}\times[\rho^2-R^2,0] \\ 0 & \textrm{if}\;\;(x,t)\in\mathbb{R}^{d+1}\setminus \mathcal{C}_R.\end{array}\right.$$
Furthermore, choose $\eta$ satisfying
$$\a{\eta_t}\lesssim \frac{1}{\rho^2}\;\;\textrm{and}\;\;\a{\nabla\eta}\lesssim\frac{1}{\rho}\;\;\textrm{on}\;\;\mathbb{R}^{d+1}.$$

Test equation (\ref{cacc_eq}) against $\eta^2(u-c)$ and use the the definition of $\eta$ and the identity
$$\nabla(\eta^2 (u-c))\cdot a\nabla u=\eta^2 \nabla u \cdot a\nabla u+2\eta (u-c) \nabla \eta\cdot a\nabla u$$
to obtain
$$\int_{\mathcal{C}_R}\eta^2\nabla u\cdot a\nabla u\lesssim \frac{1}{2}\int_{\mathcal{C}_R}(u-c)^2\partial_t\eta+\int_{\mathcal{C}_R}\eta\a{(u-c)}\a{\nabla\eta}\a{\nabla u}.$$
Therefore, following applications of H\"older's inequality and Young's inequality, and after using definition of $\eta$ and the uniform ellipticity of $a$, it follows that
$$\int_{\mathcal{C}_{R-\rho}}\a{\nabla u}^2\lesssim \frac{1}{\rho^2}\int_{\mathcal{C}_R\setminus\mathcal{C}_{R-\rho}}(u-c)^2,$$
which completes the proof.

\bibliographystyle{amsplain}
\bibliography{time-dependent}

\end{document}